\documentclass[12pt]{amsart}

\IfFileExists{srcltx.sty}{\usepackage{srcltx}}

\numberwithin{equation}{section}

\usepackage[latin1]{inputenc}
\usepackage{xspace,amssymb,amsfonts,euscript}
\usepackage{amsthm,amsmath}
\usepackage{palatino}
\usepackage{euscript}
\usepackage[dvips]{epsfig}

\usepackage{enumitem}

\usepackage[hmargin=3cm,vmargin=3.5cm]{geometry}

\usepackage{psfrag}

\input xy \xyoption {all}

\usepackage{tikz}
\usetikzlibrary{calc}
\usepackage{tikz}
\tikzset{
MyPersp/.style={scale=3,x={(-0.7cm,-0.4cm)},y={(0.8cm,-0.4cm)},z={(0cm,1cm)}},
	MyPoints/.style={fill=white,draw=black,thick}
		}

\RequirePackage{color}
\definecolor{myred}{rgb}{0.75,0,0}
\definecolor{mygreen}{rgb}{0,0.5,0}
\definecolor{myblue}{rgb}{0,0,0.65}

\RequirePackage{ifpdf}
\ifpdf
  \IfFileExists{pdfsync.sty}{\RequirePackage{pdfsync}}{}
  \RequirePackage[pdftex,
   colorlinks = true,
   urlcolor = myblue, 
   citecolor = mygreen, 
   linkcolor = myred, 
   pagebackref,
   bookmarksopen=true]{hyperref}
\else
  \RequirePackage[hypertex]{hyperref}
\fi

    \def\DM{{\mathbb{D}}}

    \def\NM{{\mathbb{N}}}
    
    \def\PM{{\mathbb{P}}}
    
    \def\RM{{\mathbb{R}}}

    \def\ZM{{\mathbb{Z}}}

  \def\ab{{\mathbf a}}  
  \def\bb{{\mathbf b}}  \def\BC{{\mathcal{B}}}
  \def\cb{{\mathbf c}}  \def\CC{{\mathcal{C}}}

    \def\MC{{\mathcal{M}}}
    
    \def\OC{{\mathcal{O}}}
    \def\PC{{\mathcal{P}}}
    \def\QC{{\mathcal{Q}}}
    \def\RC{{\mathcal{R}}}
    \def\SC{{\mathcal{S}}}

    \def\WC{{\mathcal{W}}}

    \def\ZC{{\mathcal{Z}}}

\newcommand{\nc}{\newcommand} \newcommand{\renc}{\renewcommand}

\newcommand{\rdots}{\mathinner{ \mkern1mu\raise1pt\hbox{.}
    \mkern2mu\raise4pt\hbox{.}
    \mkern2mu\raise7pt\vbox{\kern7pt\hbox{.}}\mkern1mu}}

\newcommand{\diag}{\textrm{diag}}

\def\to{\rightarrow}

\def\longto{\longrightarrow}

\nc{\triright}{\stackrel{[1]}{\to}}
\nc{\longtriright}{\stackrel{[1]}{\longto}}

\nc{\Br}{\mathcal{B}}
\nc{\HotRR}{{}_R\mathcal{K}_R}
\nc{\HotR}{\mathcal{K}_R}
\nc{\excise}[1]{}
\nc{\defect}{\text{df}}
\nc{\h}[1]{\underline{H}_{#1}}

\nc{\Ga}{\mathbb{G}_a} 
\nc{\Gm}{\mathbb{G}_m} 

\nc{\Perv}{{\mathbf{P}}}

\nc{\IH}{{\mathrm{IH}}}

\nc{\ic}{\mathbf{IC}}

\nc{\gl}{{\mathfrak{gl}}}
\renc{\sl}{{\mathfrak{sl}}}
\renc{\sp}{{\mathfrak{sp}}}

\nc{\HBM}{H^{BM}}

\DeclareMathOperator{\End}{End} 

\DeclareMathOperator{\id}{id}

\newtheorem{thm}{Theorem}[section]

\newtheorem{lemma}[thm]{Lemma}

\newtheorem{prop}[thm]{Proposition}
\newtheorem{cor}[thm]{Corollary}
\newtheorem{claim}[thm]{Claim}  
\newtheorem{conj}[thm]{Conjecture}

\theoremstyle{definition}
\newtheorem{defn}[thm]{Definition}

\newtheorem{ex}[thm]{Example}
\newtheorem{exercise}[thm]{Exercise}

\theoremstyle{remark}
\newtheorem{remark}[thm]{Remark}

\DeclareMathOperator{\Aut}{Aut}

\def\pt{{\mathrm{pt}}}

\usepackage{bbm}

\def\R{{\mathbbm R}}
\def\Z{{\mathbbm Z}}

\def\1{\mathbbm{1}}

\renewcommand{\hat}{\widehat}
\renewcommand{\tilde}{\widetilde}

\newcommand{\tX}{\tilde{X}}

\newcommand{\tOC}{\tilde{\OC}}
\newcommand{\tSC}{\tilde{\SC}}
\newcommand{\tRC}{\tilde{\RC}}

\newcommand{\tWC}{\tilde{\WC}}
\newcommand{\pa}{\partial}
\newcommand{\co}{\colon}
\newcommand{\ul}{\underline}
\newcommand{\ubr}[2]{\underbrace{#1}_{#2}}
\newcommand{\Cox}{\textrm{Cox}}
\newcommand{\tCox}{{\widetilde{\Cox}}}
\newcommand{\unor}{\textrm{un-or}}

\newcommand{\RP}{\mathbb{RP}}

\newcommand{\Int}{\mathrm{Int}}

\newcommand{\ot}{\otimes}

\newcommand{\Sal}{\textrm{Sal}}

\newcommand{\tofrom}{\rightleftarrows}

\newcommand{\ig}[2]{\vcenter{\xy (0,0)*{\includegraphics[scale=#1]{FennFig/#2}} \endxy}}
\newcommand{\igc}[2]{\begin{center} \includegraphics[scale=#1]{FennFig/#2} \end{center}}

\begin{document}

\title{Diagrammatics for Coxeter groups and their braid groups}

\author{Ben Elias} \address{Massachusetts Institute of Technology,
  Boston}

\author{Geordie Williamson} \address{Max-Planck-Institut f\"ur Mathematik, Bonn}

\begin{abstract} We give a monoidal presentation of Coxeter and braid 2-groups, in terms of decorated planar graphs. This presentation extends the Coxeter presentation. We deduce a simple
criterion for a Coxeter group or braid group to act on a category. \end{abstract}

\maketitle


\section{Introduction and preliminaries}
\label{sec-intro}

\subsection{Strictifying group actions}
\label{subsec-strictifying}

A group acting on a category is, roughly speaking, an assignment of a functor for each element of the group. This can be thought of as a categorification of the usual notion of a group
representation. The recent paradigm in categorification dictates that studying the (isomorphism classes of) functors which appear is not sufficient; one should study the natural
transformations between them. This philosophy can be found, for instance, in the seminal paper of Chuang and Rouquier \cite{CR}, which deduced strong structural results for categorified
$\mathfrak{sl}_2$ representations given the existence of a certain algebra of natural transformations.

The nuance in the work of Chuang and Rouquier was specifying an interesting algebra of natural transformations between functors. For groups, the nuance comes from the opposite goal:
showing that the algebra of natural transformations between functors corresponding to the same element of the group can be trivialized. The desired structure is a strict action of a group
on a category, where each element of the group is (compatibly) assigned a \emph{canonical} functor (see Definition \ref{defn:strictaction}). This is to be contrasted with a weak
action, where each element of the group is assigned an isomorphism class of functor. Given two words in the group which multiply to the same element, a weak action guarantees that the
corresponding compositions of functors are isomorphic, while a strict action fixes a natural transformation which realizes this isomorphism.

Here we pause to distinguish between the two most common descriptions of groups and their representations, which we call the holistic and the combinatorial. In the
holistic approach, the action of each element is given in a general way. An example is the standard representation of $GL(n)$ (or its exterior and symmetric powers), where each matrix $g
\in GL(n)$ acts via a general formula. Another example is the action of an automorphism group of a variety acting on the cohomology ring, by pullback. In the combinatorial approach, one
describes a group combinatorially by generators and relations. A representation can be defined by giving an endomorphism for each generator, and checking the relations. This can save a
great deal of labor, replacing the computation of the entire multiplication table with a manageable amount of data. Representations of Coxeter groups and their Artin braid groups are often
defined in this fashion.

These two approaches are also common when defining actions of groups on categories. The holistic approach lends itself easily to the notion of a strict action. For example, the
automorphism group of a variety acts on the derived category of coherent sheaves by pullback, and the composition of pullbacks is naturally isomorphic to the pullback of the composition.
On the other hand, given a group presentation, one could define a functor for each generator, and check an isomorphism of functors for each relation, but this would only define a weak
action. If one works directly from the definition of a strict group action, the additional data required to make this action strict is not made any simpler by the presentation: one needs
a natural transformation for each entry of the multiplication table, satisfying a host of compatibilities. This ``strictification" data is prohibitive to provide in practice, and is not
in keeping with the labor-saving combinatorial nature of the presentation.

\subsection{Braid groups and Coxeter groups}
\label{subsec-strictifyingpart2}

The primary goal of this paper is to give an explicit and efficient criterion for establishing a strict action of a Coxeter group or its braid group on a category, extending the Coxeter
presentation of said group. For instance, to make an action of the type $A$ braid group strict, one need only check a single equality: the so-called \emph{Zamolodzhikov relation}. This
improves slightly upon a similar result of Deligne \cite{Deligne} and Digne-Michel \cite{DigneMichel} for braid groups. The literature does not seem to contain any previous results on
strictifying actions of Coxeter groups.

Many topics in category theory can be more intuitively phrased using the language of topology, and strictifying a group action is a fine example. The equivalence between group
presentations and 2-dimensional cell complexes (with a single 0-cell) is well-known. Finding strictification data for this presentation is equivalent to finding a collection of 3-cells
which kill $\pi_2$ of this complex. Essentially, one is searching for a combinatorially-defined 3-skeletal approximation for the classifying space $BG$ of the group. An equivalent question
is to find an appropriate 3-skeletal approximation for its universal cover $EG$. In this paper we discuss two separate cell complexes attached to a Coxeter system, one for the Coxeter
group and the other for its associated braid group. Though the complexes are different, the corresponding strictification criteria are closely related.

To a Coxeter group $W$ one can associate a real hyperplane arrangement, and can consider the complement of these hyperplanes in the complexification $Y_W$. The
\emph{$K(\pi,1)$-conjecture}, originally due to Arnold, Brieskorn, Pham and Thom  states that
$Y_W$ should be a classifying space for the pure braid group, and thus a natural quotient $Y_W / W = Z_W$ \footnote{$Z_W$ is denoted $X_W$ in \cite{Deligne} and \cite{Sal2}. We have
reserved $X_W$ for another purpose.} should be a classifying space for the braid group. This was proven for finite Coxeter groups by Deligne \cite{Deligne}. The $K(\pi,1)$-conjecture has
also been proven for a broad class of infinite Coxeter groups. We give some further introduction in \S\ref{subsec-Kpi1}; see \cite{Paris2} for a survey.

In order to translate this topological result into a recipe for a categorical group action, one should choose a combinatorial realization of $Z_W$ as a cell complex. In his proof, Deligne
introduces a realization involving reduced expressions for elements in the Coxeter group. Correspondingly, in \cite{Deligne2}, Deligne provides an analogous criterion for a strict action
of the positive braid monoid of a finite Coxeter group. This criterion extends the \emph{positive lift presentation} of the braid group, where the generators are positive lifts of each
element in the Coxeter group. This result has been generalized to arbitrary Coxeter groups by Digne and Michel \cite{DigneMichel}, using the Garside structure on the braid group.

In \cite{Sal1}, Salvetti proved that an arbitrary hyperplane complement has the same homotopy type as a combinatorially-defined cell complex. In \cite{Sal2}, Salvetti provided an analogous
cell complex realization of the quotient $Z_W$, which he used to reprove the $K(\pi,1)$-conjecture for finite Coxeter groups. Both cell complexes are called \emph{Salvetti complexes} in
the literature; in this paper, we reserve the term for the realization of $Z_W$. Independently, Paris \cite{Paris} also used Salvetti's complexes for hyperplane complements to reprove this
result, introducing along the way a combinatorial construction for the (conjectural) universal cover of the Salvetti complex.

The Salvetti complex differs from Deligne's complex, in that the 2-skeleton corresponds to the Artin presentation of the braid group rather than the positive lift presentation. Regardless
of the validity of the $K(\pi,1)$-conjecture, the results of Digne and Michel mentioned above imply that $\pi_2(Z_W)=0$, and therefore the Salvetti complex gives a valid 3-skeletal
approximation of the classifying space. Our strictification data for the Artin presentation is the extrapolation of the 3-cells in the Salvetti complex.

In similar fashion, assuming the $K(\pi,1)$-conjecture, presumably one can use the $k$-skeleton of the Salvetti complex to concoct a strictification procedure for actions of braid groups
on $(k-2)$-categories. One can see this as a higher categorical generalization of the Coxeter presentation of a braid group; on the $k$-th categorical level, there is but a single relation
for each finite rank $k$ (standard) parabolic subgroup, a ``higher Zamolodzhikov relation." We do not pursue this any further in this paper.

We are certainly not the first to observe this type of phenomenon. Strictification data for braid group has been studied by several authors, and the Zamolodzhikov relation goes back at
least to Deligne \cite{Deligne2}. One can find the Zamolodzhikov relation in type $A$ described using Igusa pictures in Loday \cite[Figure 19]{Loday}. A similar approach for general braid
groups using the language of coherent presentations was taken in \cite{GGM}.

In contrast to the situation for braid groups, there do not seem to be cell complex realizations of the classifying spaces of Coxeter groups in the literature. However, strictification data
for Coxeter groups is essential in the authors' work on Soergel bimodules \cite{EWGR4SB}. We construct a new and somewhat unfamiliar cell complex as our 3-skeletal model for the universal
cover $EW$. We study this model by relating it to the dual Coxeter complex, in a way to be described in \S\ref{subsec-coxeterintro}.

The secondary goal of this paper is to provide diagrammatic tools for the study of Coxeter groups and braid groups. The bulk of this consists in publicizing a wonderful diagrammatic
interpretation of $\pi_2$ of a cell complex that we found in a book by Fenn \cite{Fenn} (see Remark \ref{Igusaremark}). Applying these techniques to the Salvetti complex and to our complex
for $EW$, one obtains a depiction of elements of $\pi_2$ as decorated planar graphs. This diagrammatic calculus is new for both braid groups (outside of type $A$) and Coxeter groups, and
could potentially lead to new, diagrammatic proofs of our main result.

In the remainder of this introductory chapter, we spell out the connection with topology in more detail, in order to describe our results and discuss some of the techniques we use. Then we
state our results in \S\ref{subsec-results}, using diagrammatic language, without any reference to topology. In \S\ref{subsec-applications}, we give some applications.

\subsection{Strict actions and 3-presentations}
\label{subsec-3present}

\begin{defn} For a group $G$, let $\Omega G$ be the monoidal category defined as follows. The objects consist of the set $G$, and the only morphisms are identity maps $\id_g$ for each $g
\in G$. The monoidal structure on objects is given by the group structure on $G$, and the monoidal structure on morphisms is uniquely determined. \end{defn}

\begin{defn} Given a category $\CC$, let $\Aut(\CC)$ denote the monoidal category whose objects are autoequivalences of $\CC$, and whose morphisms are invertible natural transformations.
The monoidal structure is given by composition of functors. \end{defn}

\begin{defn} \label{defn:strictaction} A \emph{strict group action} of $G$ on a category $\CC$ is a monoidal functor $\Omega G \to \Aut(\CC)$. \end{defn}

\begin{remark} \label{remark:usualdefinitionstrict} The usual definition of a strict group action involves providing a functor $F_g$ for each $g \in G$, isomorphisms $a_{g,h} \co F_g
\circ F_h \to F_{gh}$, and an isomorphism $\epsilon \co F_e \to \1_{\CC}$ for the identity element $e \in G$, satisfying some natural compatibilities, including an associativity
compatibility. The isomorphisms $a_{g,h}$ are the image of the unique morphism $g \ot h \to gh$ in $\Omega G$, and the isomorphism $\epsilon$ is the fixed isomorphism between monoidal
identities given as part of the data of a monoidal functor. \end{remark}


Suppose that $G$ has a presentation $\PC = (\SC,\RC)$ with generators $\SC$ and relations $\RC$. As discussed above, it is common in the literature to define a weak group action by giving
an invertible functor $F_s$ for each $s \in \SC$, and checking an isomorphism for each $r \in \RC$. Doing this is implicitly defining a monoidal functor from the monoidal category $\Omega
\PC$ defined below. The inverse functor $F_s^{-1}$ is not usually given explicitly, but one can choose $F_s^{-1}$ to be any inverse to $F_s$.

\begin{defn} For a presentation $\PC = (\SC,\RC)$, let $\Omega \PC$ be the monoidal category defined as follows. Its objects are words in the letters $\SC \cup \SC^{-1}$, with monoidal
structure given by concatenation (we let $\1$ denote the monoidal identity, the empty word). Its morphisms are monoidally generated by the following maps: \begin{itemize} \item (Cups and
Caps) Inverse isomorphisms $ss^{-1} \tofrom \1$ for each $s \in \SC$, as well as inverse isomorphisms $s^{-1}s \tofrom \1$. \item (Relation symbols) Inverse isomorphisms $r \tofrom \1$ for
each word $r \in \RC$. \end{itemize} One imposes the following relations: \begin{itemize} \item The generating ``inverse isomorphisms" are actually inverse isomorphisms. \item Cups and
caps form the units and counits of adjunction between the biadjoint functors $s \ot (\cdot)$ and $s^{-1} \ot (\cdot)$. \item The relation symbols are cyclic with respect to these
biadjunction structures. \end{itemize} \end{defn}

\begin{remark} This remark is for the reader unfamiliar with biadjunction and cyclicity. For a broader introduction to biadjunction, cyclicity, and the associated diagrammatic notation,
see \cite{LauSL2,LauDiagrams}.

Observe that any invertible functor $F \in \Aut(\CC)$ is both left and right adjoint to its inverse functor $F^{-1}$. Adjunction is a structure, not a property, and a \emph{biadjunction}
is a choice of both a left and right adjunction between $F$ and $F^{-1}$. Having chosen inverse isomorphisms $F_s F_s^{-1} \tofrom \1$, there is a unique choice for the isomorphisms
$F_s^{-1} F_s \tofrom \1$ such that the isomorphisms also provide (the units and counits of) a biadjunction. Biadjunctions occur more generally than between inverse functors, but this
will suffice for our purposes.

\emph{Cyclicity} is a property of a general morphism (i.e. natural transformation) between compositions of functors equipped with biadjunctions, stating that this morphism is somehow
compatible with the right versus the left adjunction. An explicit statement of this compatibility can be found (in diagrammatic language) later in this paper. Identity morphisms are
axiomatically cyclic, but general morphisms need not be cyclic.

In a monoidal category (like $\Aut(\CC)$), the biadjoint of an object, if it exists, is well-defined up to unique isomorphism, but nonetheless this assignment of a biadjoint to each
object need not be functorial. In a \emph{pivotal} category, there exists a duality functor $\DM$ sending each object to a biadjoint, and equipped with a natural isomorphism $\phi \co \1
\to \DM^2$. One can always adjust the functor $\DM$ up to isomorphism to guarantee that $\DM^2$ and $\1$ agree on objects, but it need not be the case that $\phi$ is the identity map. If
$\phi$ is the identity map, the category is \emph{strictly pivotal}. Equivalently, every morphism is cyclic, so the category is also called \emph{cyclic biadjoint}. This happens
frequently in geometric and algebraic examples of group actions. It also happens in fundamental 2-groupoids, the topological framework of this paper. \end{remark}

\begin{remark} Note that a functor $\Omega \PC \to \Aut(\CC)$ is not quite as general as a weak group action defined by generators and relations, because there need not exist relation
symbols which are cyclic. Any weak action which can be extended to a strict action will certainly satisfy cyclicity. \end{remark}

\begin{remark} What is the difference between a categorical action of a monoid where the generators happen to act by invertible functors, and a categorical action of its associated group?
From the definition of a weak categorical action, there is no difference. Philosophically, however, one might desire some new condition which connects the new inversion structure with the
existing relation isomorphisms. Said another way, one now has a host of new relations obtained by conjugating existing relations, and one might expect these to be somehow mutually
compatible. We believe that cyclicity is precisely the correct structure one should impose. \end{remark}

Every morphism in $\Omega \PC$ is an isomorphism. The isomorphism classes of objects in $\Omega \PC$ can be identified with $G$, and there is a monoidal functor $\Omega \PC \to \Omega G$.
However, endomorphism spaces in $\Omega \PC$ can be quite large, so this functor is not faithful. Because of biadjunction, every endomorphism space is a principal space for the group
$\End(\1)$.

\begin{defn} Let $\ZC$ be a chosen subset of $\End(\1)$ within the category $\Omega (\SC,\RC)$, for some group presentation $(\SC,\RC)$ of $G$. We call $\PC=(\SC,\RC,\ZC)$ a
\emph{3-presentation} of $G$, and we simply write $(\SC,\RC)$ when $\ZC$ is empty. We let $\Omega \PC$ denote the quotient of $\Omega (\SC,\RC)$ by the relations $z = \id_{\1}$ for each $z
\in \ZC$. \end{defn}

For any 3-presentation $\PC$ of $G$, there is still a monoidal functor $\Omega \PC \to \Omega G$. When $\ZC$ generates the group $\End(\1) \subset \Omega (\SC,\RC)$, then morphism spaces
in $\Omega \PC$ are trivial, consisting only of identity maps, and the functor to $\Omega G$ is an equivalence. We call such a 3-presentation \emph{acyclic}. In the acyclic case,
giving a monoidal functor $\Omega \PC \to \Aut(\CC)$ is equivalent to giving a strict action of $G$, but it has a different recipe: provide a functor $F_s$ and its biadjoint inverse for
each $s \in \SC$, provide an isomorphism and its inverse for each $r \in \RC$, and check a relation for each $z \in \ZC$.

This recipe need not be interesting or useful. We do not expect there to be a general method to extend an arbitrary $2$-presentation of a group into an acyclic $3$-presentation in a
\emph{useful} way.

\begin{ex} Every group has a universal presentation where $\SC = G$ and $\RC$ consists of relations stating that $g \cdot h = (gh)$. The corresponding monoidal category $\Omega
(\SC,\RC)$ has isomorphisms $a_{g,h}$ as in Remark \ref{remark:usualdefinitionstrict}, but no compatibility requirements. Letting $\ZC$ be the set of associativity requirements (one for
each triple $g,h,k \in G$), one has that $(\SC,\RC,\ZC)$ is acyclic. This is the universal 3-presentation of a group, and one could say that it is the only uninteresting
3-presentation of the group, since it does not reduce the labor required to construct a strict action. \end{ex}

\begin{ex} Given any presentation $(\SC,\RC)$, the 3-presentation $(\SC,\RC,\End(\1))$ is its universal acyclic extension. This example is also not very interesting or useful,
because computing $\End(\1)$ and checking a relation for each element of $\End(\1)$ can be prohibitive. \end{ex}

In this paper we will give \emph{interesting} examples of acyclic 3-presentations, extending the usual presentations of Coxeter groups and their braid groups. Said another way, we find
an interesting presentation of the uninteresting monoidal category $\Omega G$ for these groups.

\subsection{Diagrammatics and Topology}
\label{subsec-topology}

To any topological space $X$ one may associate its fundamental 2-groupoid $\pi(X)_{\le 2}$ (see e.g. \cite[8.2]{BaezLauda}) \footnote{Technically, this is called the ``homotopy bigroupoid"
in \cite{BaezLauda}.}. In this 2-category, the objects are the points of $X$, the 1-morphisms from $x$ to $y$ are given by paths from $x$ to $y$, and the 2-morphisms are given by ``paths
of paths'' up to homotopy. Our notation is intended to suggest that $\pi(X)_{\le 2}$ is a truncation of the fundamental $\infty$-groupoid $\pi(X)$, which encodes the homotopy type of $X$. When $X$ is a cell complex, $\pi(X)_{\le 2}$ only depends on the 3-skeleton $X^3 \subset X$.

There is an explicit diagrammatic interpretation of $\pi(X)_{\le 2}$ for a cell complex $X$, known in the literature as \emph{Igusa's pictures} \cite{Igusa} (see Remark \ref{Igusaremark}).
After fixing some additional data, one can encode the structure of the cell complex combinatorially. A sufficiently nice map $D^2 \to X^2$ is depicted as a decorated oriented planar graph;
we call such a map \emph{pictorial}. Any map $D^2 \to X$ is homotopic to a pictorial map. Moreover, there is a list of relations which account for all homotopies between pictorial maps,
essentially arising from Morse theory. This diagrammatic calculus can be viewed as a tool that takes a 3-skeleton of a cell complex $X$ and returns a combinatorially-defined 2-category
$\Pi(X)_{\le 2}$, which is equivalent to $\pi(X)_{\le 2}$. The objects of $\Pi(X)_{\le 2}$ are the 0-cells of $X$; the generating 1-morphisms are the 1-cells and their formal (biadjoint)
inverses; the generating 2-morphisms are the 2-cells and their formal inverses, along with units and counits of biadjunction; and the relations between 2-morphisms are given by the
3-cells, along with some general relations (inverses are inverses, other 2-morphisms are cyclic). A 2-morphism in $\Pi(X)_{\le 2}$ will be represented by a decorated planar graph, whose
regions are labelled by 0-cells, whose edges are labelled by 1-cells with an orientation, and whose vertices are labelled by 2-cells with some additional data.

Let $\PC=(\SC,\RC)$ be a presentation of $G$. In a standard way, this is also a recipe for a 2-complex $X_\PC$ with a single 0-cell, for which $\pi_1(X_\PC) \cong G$. The corresponding
monoidal category $\Pi(X_\PC)_{\le 2}$ is the category $\Omega \PC$ defined above. Similarly, for a 3-presentation $\PC = (\SC,\RC,\ZC)$ there is a 3-complex $X_\PC$,
and $\Pi(X_\PC)_{\le 2}$ equals $\Omega \PC$. The 3-presentation is acyclic if and only if $\pi_2(X_\PC)$ is trivial, in which case $\Pi(X_\PC)_{\le 2} \cong \Omega G$, and $X_\PC$ is the
3-skeleton of some realization of the classifying space $BG = K(G,1)$.

The presentation of a Coxeter group has a number of symmetries (though the presentation of its braid group does not). Exploiting these symmetries, we can modify this construction
of $\Pi(X_W)_{\le 2}$ to produce a simpler diagrammatic calculus. For example, when constructing a 2-morphism as a decorated planar graph, one need not specify the orientation on the
edges; think of this as using the relation $s^2=1$ to canonically identify $s$ and $s^{-1}$.

\subsection{The $K(\pi,1)$ conjecturette}
\label{subsec-Kpi1intro}

Let $(W,\SC)$ be a Coxeter system and $B_W$ be the corresponding Artin braid group.

In chapter \S\ref{sec-braiddiag} we define a 3-presentation $\PC_{B_W}$ of $B_W$. The corresponding 3-complex $X_{\PC_{B_W}}$ is the 3-skeleton of a cell complex we shall call $X_{B_W}$,
which is the Salvetti complex. Further discussion of this complex can be found in \S\ref{subsec-Kpi1}. As discussed earlier in the introduction, $X_{B_W}$ is the subject of a famous
conjecture.

\begin{conj} (The $K(\pi,1)$-conjecture) $X_{B_W}$ is the classifying space of $B_W$. \end{conj}

However, the fact that $\PC_{B_W}$ is acyclic is equivalent to a much weaker condition, which we call the \emph{$K(\pi,1)$-conjecturette}.

\begin{prop} (The $K(\pi,1)$-conjecturette) $\pi_2(X_{B_W})$ is trivial. \end{prop}

Via the work of Salvetti \cite{Sal1}, this becomes a question about hyperplane complements, and Digne-Michel's generalization \cite[\S 6]{DigneMichel} of Deligne's finite-type arguments
gives a proof of the $K(\pi,1)$-conjecturette for all Coxeter groups. We will quote this result henceforth. However, we believe our diagrammatic tools should allow for an elementary and
direct proof, which unfortunately has not yet materialized.

\subsection{Modified Coxeter complexes and half-skeletons}
\label{subsec-coxeterintro}

In \S\ref{sec-coxdiag} we define a 3-presentation $\PC_W$ of $W$. The corresponding 3-complex $X_W$, which we call (the 3-skeleton of) the \emph{modified Coxeter complex}, is not a
familiar topological space. However, a natural $W$-fold cover $\tX_W$ of $X_W$ has an equivalent fundamental $2$-groupoid to the 3-skeleton of the completed dual Coxeter complex $\tCox_W$.
The completed dual Coxeter complex will be discussed further in \S\ref{sec-coxbackground}. Thankfully, $\tCox_W$ is known to be contractible, therefore giving the proof that $\pi_2(\tX_W)
\cong \pi_2(X_W)$ is trivial.

The classifying space $BW$ (for which $X_W$ is supposed to be a model) is a quotient of its universal cover $EW$ (for which $\tX_W$ is supposed to be a model). The space $EW$ must satisfy
two conditions: it must be contractible, and it must admit a free $W$-action. When trying to build $EW$ as a cell complex, one is torn between these two goals. Perhaps the neatest approach
is to alternate between them, first constructing a contractible space, then extending it until it admits a free $W$-action, then extending it to make it contractible again, and so forth.
This leads to the concept of \emph{half-skeletons}, which we use to prove the result about $\tX_W$ and $\tCox_W$. Half-skeletons are not intended to be a complete theory, just a heuristic
organizational tool.

Let us illustrate the approach using the simplest example of a Coxeter group, $W = \ZM/2\ZM$, an example which is treated in more detail within the body of the paper. The classifying space
of $W$ is $\RM \PM^\infty$, with universal cover $S^\infty$. The standard cell complex construction of $S^\infty$ has two $k$-cells for each $k \ge 0$, such that the
$k$-skeleton is $S^k$. Note that $S^k$ admits a free $W$-action compatible with the cell decomposition, but is not contractible.

Now take the $k$-skeleton $S^k$ and attach just one of the two $(k+1)$-cells; one obtains the disk $D^{k+1}$, which is contractible, but does not admit a free $W$-action. This is the $(k +
0.5)$-skeleton of $EW$. To get from the $(k + 0.5)$-skeleton to the $(k+1)$-skeleton, one attaches an additional $(k+1)$-cell, but along an attaching map which is nulhomotopic in the $(k +
0.5)$-skeleton (unsurprisingly). Topologically, this operation is just wedging with $S^{k+1}$, and hence does not change $\pi_l$ of the space for any $l \le k$.

Let $X^k$ or $X^{k + 0.5}$ denote such a skeleton, for $k \in \NM$. To reiterate, this setup is designed so that: \begin{itemize} \item $X^k$ admits a free cellular $W$-action. \item $X^{k+0.5}$ is
contractible (or at least has trivial fundamental $k$-groupoid). \item To get from $X^{k+0.5}$ to $X^{k+1}$, one attaches a set of $(k+1)$-cells (along attaching maps which are necessarily
nulhomotopic). \end{itemize} In particular, this guarantees that $X^{k+0.5}$ and $X^{k+1}$ have equivalent fundamental $k$-groupoids, so that $\pi_l(X^{k+1})=0$ for $l \le k$.

In \S\ref{sec-coxdiag}, for a general Coxeter group $W$, we construct a $2.5$-skeleton and a $3$-skeleton for $EW$. The $2.5$-skeleton will naturally deformation retract to the completed
dual Coxeter complex $\tCox_W$, and is therefore contractible. The $3$-skeleton is exactly $\tX_{\PC_W}$ for our chosen 3-presentation. This explains why $\tX_{\PC_W}$ and $\tCox_W$ have
equivalent fundamental $2$-groupoids $\pi_{\le 2}$.

We do not propose a combinatorial method to construct higher skeletons and half-skeletons for $EW$, largely because the ``diagrammatic" technology for understanding higher fundamental
groupoids is undeveloped.

\subsection{Results}
\label{subsec-results}

Let $(W,S)$ be a Coxeter system and let $B_W$ be the corresponding Artin braid group. We now describe our presentations of $\Omega B_W$ and $\Omega W$ without any reference to topology.

\begin{defn} \label{defn:BCdiag}
Let $\BC_{\diag}$ be the monoidal category with
\begin{enumerate}
\item objects -- words in $S \cup S^{-1}$, or equivalently, sequences of
  oriented dots on a line colored by $S$. \igc{1.5}{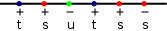}
\item morphisms -- planar strip diagrams, generated by oriented
  cups, caps and two-colored $2m$-valent vertices. (In the example below, $m_{st}=3$ and $m_{su}=2$. We will continue to use these to exemplify the general case. There is no generator when $m_{st} = \infty$.) \igc{1.5}{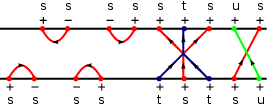} \end{enumerate}
Morphisms are taken modulo the relations below. Each relation holds for any valid ``coloring," i.e. any valid labeling of the edges by elements of $S$.
\begin{enumerate}[resume]
\item 
The \emph{standard relations}:
\begin{subequations} \label{Fennrelationsintro}
\begin{equation} \ig{1.75}{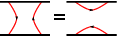} \end{equation}
\begin{equation} \ig{1.5}{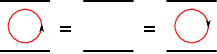} \end{equation}
\begin{equation} \ig{1.5}{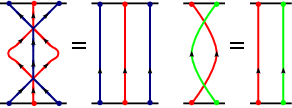} \end{equation}
\end{subequations}
\item  The \emph{isotopy relations}:
	\begin{subequations} \label{isotopyrelationsintro}
	\begin{equation} \ig{1}{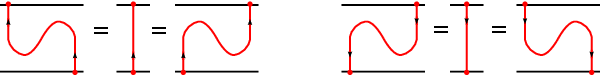} \label{biadjointintro} \end{equation}
	\begin{equation} \ig{1}{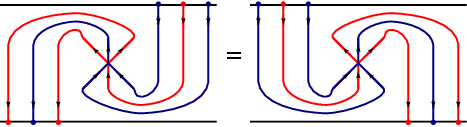} \label{2mcyclicintro} \end{equation}
(This picture illustrates the case $m = 3$. We require a similar relation for any $2m$-valent vertex.)
	\end{subequations}
\item
The \emph{generalized Zamolodzhikov relations}, one for each finite (standard) parabolic subgroup of rank 3:
\begin{subequations} \label{Zamrelationsintro}
\begin{equation} \textrm{Type}\;A_1 \times I_2(m): \quad \quad \ig{1}{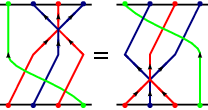} \end{equation}
\begin{equation} \textrm{Type}\;A_3: \quad \quad \ig{1}{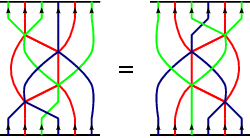} \label{A3intro} \end{equation}
\begin{equation} \textrm{Type}\;B_3: \quad \quad \ig{1}{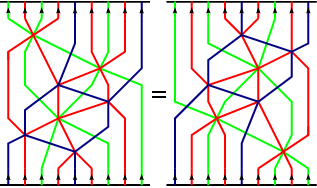} \end{equation}
\begin{equation} \textrm{Type}\;H_3: \quad \quad \ig{1}{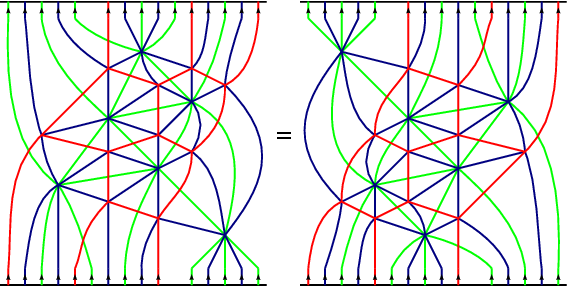} \end{equation}
\end{subequations}
\end{enumerate}

This ends the definition. \end{defn}

\begin{remark} The isotopy relations are equivalent to the statement that one can consider these embedded planar graphs up to isotopy. \end{remark}

\begin{remark} In type $A_3$, the relation \eqref{A3intro} is related to the Zamolodzhikov relation. For this reason, we call the relations in \eqref{Zamrelationsintro} generalized
Zamolodzhikov relations. \end{remark}

\begin{thm} \label{mainthmbraid} The obvious monoidal functor $\BC_\diag \to \Omega B_W$ is an equivalence of categories. \end{thm}

\begin{cor} To define a strict action of $B_W$ on a category $\CC$ is equivalent to giving the following data: \begin{enumerate}
\item Functors $F_s$ and $F_s^{-1}$ for each $s \in S$.
\item Natural transformations $F_s F_s^{-1} \tofrom \1$ and $F_s^{-1} F_s \tofrom \1$ for each $s \in S$. 
\item For each $s,t \in S$ with $m_{st}$ finite, natural transformations $F_s F_t F_s \ldots  \tofrom F_t F_s F_t \ldots$; here each expression has length $m_{st}$.
\end{enumerate}
This data is subject to the following conditions: \begin{enumerate}[resume]
\item The pairs of natural transformations above are inverse isomorphisms.
\item The natural transformations identifying $F_s^{-1}$ as the inverse of $F_s$ form a biadjoint structure.
\item The natural transformations $F_s F_t F_s \ldots  \tofrom F_t F_s F_t \ldots$ are cyclic with respect to this biadjoint structure.
\item The generalized Zamolodzhikov relations hold. \end{enumerate}
These conditions correspond to \eqref{Fennrelationsintro}, \eqref{biadjointintro}, \eqref{2mcyclicintro}, and \eqref{Zamrelationsintro} respectively.
\end{cor}

Note that the generalized Zamolodzhikov relations do not involve the functors $F_s^{-1}$ at all, and can be checked without needing to fix these inverse functors and the cups and caps.
Given a collection of invertible functors $F_s$ acting on a suitably nice category (i.e. one with functorial biadjoints) and satisfying the generalized Zamolodzhikov relations, one can
cook up the rest of the data.

\begin{defn} Let $\WC_\diag$ be the category defined as in Definition \ref{defn:BCdiag}, except without any orientations. In other words, objects are words in $S$, or equivalently,
sequences of (unoriented) colored dots on a line. Morphisms are diagrams up to isotopy, generated by (unoriented) cups, caps, and $2m$-valent vertices, modulo the unoriented versions of
all the relations above. \end{defn}

\begin{thm} \label{mainthmcox} The obvious monoidal functor $\WC_\diag \to \Omega W$ is an equivalence of categories. \end{thm}

There is a monoidal functor $\BC_\diag \to \WC_\diag$, which on objects sends both $s^{+}$ and $s^{-}$ to $s$, and on morphisms sends an oriented diagram to its unoriented version. This
categorifies the quotient map $B_W \to W$.

\subsection{Organization of the paper}
\label{subsec-organization}

We have divided the paper into two parts: the purely topological, and the Coxeter-theoretic.

The first half of this paper will give an exposition of the diagrammatic approach to $\pi_{\le 2}$ (\S\ref{sec-Fenn}), and of the modifications one can perform in the presence of
symmetries in a group presentation (\S\ref{sec-modifiedFenn}). It is independent of the rest of the paper, although some of the key examples are motivated by the Coxeter complex (see
\S\ref{sec-coxbackground}). The modified construction uses in some sense the concept of half-skeletons, as introduced above.

The second half of the paper begins (\S\ref{sec-coxbackground}) by providing background on Coxeter groups, Coxeter complexes, and Salvetti complexes. In \S\ref{sec-braiddiag} we define a
3-presentation of the braid group, whose diagrams agree with $\BC_{\diag}$ above, and whose 3-complex agrees with the 3-skeleton of a quotient of the Salvetti complex. The $K(\pi,1)$-conjecturette states that $\pi_2$ of the Salvetti complex is trivial, which implies Theorem \ref{mainthmbraid}. This chapter only requires \S\ref{sec-Fenn}, not needing the
modified diagrams from \S\ref{sec-modifiedFenn}. In \S\ref{sec-coxdiag} we define a 3-presentation of the Coxeter group, whose diagrams agree with $\WC_\diag$ above. We show that its
3-complex has a 2.5-skeleton which deformation retracts to the dual Coxeter complex, as discussed in \S\ref{subsec-coxeterintro}. This proves Theorem \ref{mainthmcox}. This chapter will
require the modified diagrams from \S\ref{sec-modifiedFenn}.

The reader interested only in the braid group can safely skip \S\ref{sec-modifiedFenn} and \S\ref{sec-coxdiag}, and can ignore any mention of half-skeletons. The reader interested only in
the Coxeter group can safely skip \S\ref{sec-braiddiag}.

\subsection{Applications and further directions}
\label{subsec-applications}

The authors came to this topic in their study of Soergel bimodules \cite{Soe5}, which provide a categorical action of the Hecke algebra of $W$. Certain complexes of Soergel bimodules
(known as \emph{Rouquier complexes}, see \cite{RouquierBraid}) give braid group actions that, after localization, become Coxeter group actions. The description of strict Coxeter group
actions given in this paper also gives a presentation of the monoidal category of localized Soergel bimodules, which is essential to our description of the category of Soergel bimodules
\cite{EWGR4SB}.

Braid group actions on categories appear to be ubiquitous in modern geometric representation theory. These braid group actions are defined using the Coxeter presentation, and so the
authors have typically (understandably) neglected to make these actions strict. However, many of these examples have since been proven to be strict. Examples of such braid group actions
include: Bondal-Kapranov's construction of mutations on triangulated categories \cite{BondalKapranov}; Brou\'e-Michel's construction in Deligne-Lusztig theory \cite{BroueMichel}; Khovanov's homology of tangles \cite{KTang};
Seidel-Thomas twists around (-2)-curves on derived categories of coherent sheaves \cite{ST} and generalizations \cite{BezRiche}; and braid group actions via shuffling functors in highest
weight representation theory. Deligne had the Bondal-Kapranov and Brou\'e-Michel constructions in mind as applications when he gave his criterion for a braid group action \cite{Deligne2}.
Indeed, Deligne's criterion does seem to be sufficient for many constructions.

Other categorical actions of the braid group, such as those arising in categorical actions of Kac-Moody algebras, have not yet been proven to be strict.

Occasionally one can show that a certain space of natural transformations is only one-dimensional, and can conclude that strictification data exists without being forced to provide it
explicitly. This is the approach taken by Rouquier \cite{RouquierBraid} and Khovanov-Thomas \cite{KT}. In these cases, it is now trivial to make the strictification data explicit. For
Rouquier complexes in type $A$, an explicit approach was taken in \cite{EKr}.

Many categorical actions of the braid group descend to actions of the Hecke algebra on the Grothendieck group, such as those arising in highest weight representation theory. One expects
these to arise from a categorical action of the Hecke algebra, which would then imply the strictness of these actions. Spherical twists such as those in \cite{ST} can be investigated using
technology developed by Joseph Grant \cite{GrantPeriodic, GrantLifts}. In upcoming work of the first author and Grant, we will demonstrate the connection between certain actions by
spherical twist and categorical Hecke algebra actions.

In type $A$, words in the braid group are themselves topological objects, and one has the notion of braid cobordisms between such words. Braid cobordisms also have a description by
generators and relations due to Carter-Saito \cite{CarterSaito}: the morphisms are called \emph{movies}, and the relations \emph{movie moves}. Not all movies are invertible, however.

\begin{cor} To give a strict braid group action in type $A$ is the same data as an action of the invertible braid cobordism category. \end{cor}

\begin{proof} In fact, our description of $\BC_\diag$ by generators and relations agrees with that of Carter-Saito for the invertible part of their braid cobordism category. Movie moves 3,
5, 6, and 7 correspond to the standard relations \eqref{Fennrelationsintro}; movie moves 1, 2, and 8 correspond to the isotopy relations \eqref{isotopyrelationsintro}; and movie moves 4, 9,
and 10 correspond to the Zamolodzhikov relations \eqref{Zamrelationsintro}. \end{proof}

In \cite{EKr}, the first author and Daniel Krasner proved that the entire braid cobordism category acts on Rouquier complexes, not just the invertible part. Which other actions
admit such an extension, and what corresponds to the non-invertible braid cobordisms in other types, are both interesting questions.

Finally, one should note that Theorems \ref{mainthmbraid} and \ref{mainthmcox}, which are stated in purely diagrammatic language, could admit purely diagrammatic proofs. This can be
accomplished in a variety of special cases (e.g. Coxeter groups in type $A$, by an argument similar to the one used in \cite{EKh}). However, no general diagrammatic proof currently exists.

Affine Weyl groups are Coxeter groups, but also admit another presentation, the \emph{loop presentation}. Often, weak categorical actions of affine Weyl groups are given using the loop
presentation rather than the Coxeter presentation, and thus different strictification data are required. We do not consider this (interesting) question in this paper.

\emph{Acknowledgements:} We would like to thank Ruth Charney and Jean Michel for useful correspondence. The first author was supported by NSF Postdoctoral Fellowship DMS-1103862.

\part{Topology and diagrammatics}
\label{pt:top}

\section{Igusa Diagrams}
\label{sec-Fenn}

The goal of this chapter will be to take a cell complex $X$ and construct, by generators and relations, a cyclic biadjoint $2$-category $\Pi(X)_{\le 2}$ which is equivalent to the
fundamental 2-groupoid $\pi(X)_{\le 2}$ of $X$. We assume the reader is familiar with diagrammatic interpretations of cyclic biadjoint $2$-categories, for which an excellent introduction
can be found in \cite{LauSL2,LauDiagrams}.

The $2$-morphisms in $\Pi(X)_{\le 2}$ should be combinatorial encodings of maps $D^2 \to X$. Note that $\Pi(X)_{\le 2}$ only depends on the 3-skeleton of $X$, so we may assume that $X$ is
a 3-complex. We follow the procedure described by Roger Fenn in his book \cite{Fenn}. First, enrich the notion of a 3-complex to make it more combinatorial, by adding a small amount of
data for each cell and placing minor restrictions on attaching maps, none of which is significant up to homotopy equivalence. Given an enriched 2-complex, certain
planar diagrams, \emph{Igusa diagrams}, can be used to encode (nice) maps $D^2 \to X^2$, such as the attaching maps of the 3-cells. Finally, one lists the relations
between diagrams which correspond to homotopy in an enriched 3-complex.

To any $3$-presentation $\PC = (\SC,\RC,\ZC)$ of a group $G$, one can associate a 3-complex $X_\PC$ for which $\pi_1(X_\PC)=G$. One can also construct the universal cover $\tX_{\PC} \to
X_{\PC}$ in such a way that the action of $G$ is inherent from the cell complex structure on $\tX_{\PC}$. We discuss the diagrammatics for the corresponding fundamental 2-groupoids below.

Fenn's exposition is highly recommended. We give a quick summary, following sections 1.2, 2.3 and 2.4 of \cite{Fenn}. Fenn's discussion requires that $X$ have a unique 0-cell, but it is
straightforward to generalize to a cell complex with multiple 0-cells, as we do below. It is also straightforward to reorganize everything into a 2-category, with one object for each
0-cell.

\begin{remark} \label{Igusaremark} This diagrammatic interpretation of $\pi(X)_{\le 2}$ for a cell complex $X$ is credited to Whitehead by Igusa \cite[Remark following Proposition
7.4]{Igusa}. Subsequent papers (e.g. \cite{Loday,Wagoner,KapSai}) call these diagrams ``Igusa pictures."  It seems likely that Fenn independently discovered this diagrammatic description \cite{Fenn}. \end{remark}

\subsection{Cell complexes and pictorial maps}
\label{subsec-cellcomplexes}

\begin{defn} An \emph{(enriched) 3-complex} will be the following data. \begin{itemize}
	\item A set of 0-cells $\OC$.
	\item A set of 1-cells $\SC$, viewed as oriented edges $D^1$ between 0-cells. For each $s \in \SC$, we fix a point $\hat{s} \in \Int(s)$, and let $\hat{\SC}=\{\hat{s}\}_{s \in \SC}$.
	\item A set of 2-cells $\RC$, viewed as oriented disks $D^2$ attached along their boundary to the above oriented graph. We assume each attaching map is \emph{pictorial}, in a sense to be defined shortly. For each $r \in \RC$, we fix a point $\hat{r} \in \Int(r)$ and a point $\pt_r \in \pa(r)$, and let $\hat{\RC} = \{\hat{r}\}_{r \in \RC}$.
	\item A set of 3-cells $\ZC$, viewed as oriented balls $D^3$ attached along their boundary to the above 2-skeleton. We assume each attaching map is \emph{pictorial}, in a sense to be defined shortly. For each $z \in \ZC$ we fix a point $\pt_z \in \pa(z)$. \end{itemize} \end{defn}

We now define diagrammatic, combinatorial ways to encode maps from $D^1 \to X^1$ and $D^2 \to X^2$.

\begin{defn} A \emph{line diagram} is an interval $D^1$ decorated as follows. A finite number of points in $\Int(D^1)$ are labeled with an element of $\SC$ and an orientation $\pm$, i.e.
with an element of $\SC \cup \SC^{-1}$. In this paper we associate a color to each $s \in \SC$, and we refer to this labeling as a ``coloring." The regions between those points are
labeled with elements of $\OC$. The region to the left of a point labelled $s^{+}$ (resp. $s^{-}$) must be the source (resp. target) of the oriented edge $s$, and the region to the right
must be the target (resp. source).\end{defn}

To a line diagram $f$ we have a word $w(f)$ in the letters $\SC \cup \SC^{-1}$, which determines the line diagram uniquely. We also have a word $o(f)$ in the letters $\OC$. Clearly
$w(f)$ determines $o(f)$, while $o(f)$ determines $w(f)$ so long as $X^1$ has no loops or double edges.

\begin{ex} This is an example where the $0$-cells are labelled $\{\ab,\bb,\cb\}$ and the $1$-cells are labelled $\{r,g,b\}$ for \emph{r}ed, \emph{g}reen and \emph{b}lue. We will continue
this example below. We have drawn a line diagram whose word $w$ is $brg^{-1}brr^{-1}$. \igc{2}{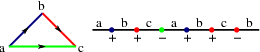} \end{ex}

Suppose that $f \co D^1 \to X^1$ is a map. It is \emph{represented} by a given line diagram if \begin{itemize} \item For each $s \in \SC$, $f^{-1}(\hat{s})$ is the collection of points
colored $s$. \item Each point colored $s$ in $D^1$ has a neighborhood which maps homeomorphically to a neighborhood of $\hat{s}$. The sign on that point is $+$ if the homeomorphism
preserves orientation, and $-$ if it reverses it. \end{itemize} Note that each connected component of $X^1 \setminus \hat{\SC}$ is star-shaped and deformation retracts to a single 0-cell
$o \in \OC$. The conditions above imply that the entire region labelled $o$ in $D^1$ will map to the corresponding connected component. Also note that the endpoints $\pa(D^1)$ can not map
to $\hat{\SC}$.

\begin{defn} A map $D^1 \to X^1$ is \emph{pictorial} if it is represented by some line diagram. A map $(S^1,\pt) \to X^1$ is \emph{pictorial} if the corresponding map $D^1 \to X^1$ is
pictorial, given by identifying $\pa(D^1)$ with $\pt \in S^1$. \end{defn}

It is easy to modify the notion of a line diagram to obtain that of a \emph{circle diagram}, representing a map $(S^1,\pt) \to X^1$. We keep track of the marked point with a tag. One can
\emph{flip} a line or circle diagram, which will invert all the orientations, and will correspond to the obvious precomposition with the flip map $D^1 \to D^1$ or $S^1 \to S^1$.

\begin{ex} This is a loop with word $g^{-1}brr^{-1}r$ based at $\cb$, and its flip $r^{-1}rr^{-1}b^{-1}g$. \igc{1.5}{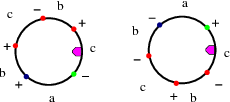} \end{ex}

Any line diagram is clearly realized by some map $D^1 \to X^1$. Any two pictorial maps $D^1 \to X^1$ with isotopic line diagrams are clearly homotopic (via a homotopy sending $\pa(D^1) \to
X^1 \setminus \hat{\SC}$). In the definition of an enriched 3-complex, the attaching map of a 2-cell $r$ is assumed to be pictorial, and thus has a circle diagram; the marked point $\pt_r \in \pa(r)$ corresponds to the tag.

\begin{defn} A \emph{disk diagram} is a particular kind of oriented planar graph in the disk $D^2$. Each edge of the graph is colored with some $s \in \SC$, and each region is labeled with
some $o \in \OC$, compatible with the orientations on edges via a ``left-handed rule." Edges may run to the boundary $\pa(D^2)$, yielding a circle diagram on the boundary (see example for
orientation rules). Edges need not meet any vertices, forming circles, or arcs at the boundary. Each vertex is labelled with an element of $\RC$ and an orientation $\pm$, i.e. with an
element of $\RC \cup \RC^{-1}$. A small circle around a vertex labelled $r^{+}$ (resp. $r^-$) must yield the circle diagram of $r$ (resp. the flip of the circle diagram of $r$). A
\emph{disk diagram with marked points} on the boundary is exactly that, with the additional assumption that the marked points do not meet the edges. \end{defn}

\begin{ex} To the $1$-skeleton of the previous examples we have glued a $2$-cell $w$ along $b^{-1}brr^{-1}$ (based at the $0$-cell $\bb$) and another $v$ along $gr^{-1}b^{-1}$ (based at
$\ab$). Now we have constructed a map from the disk which uses $w$ and $v^{-1}$. \igc{2}{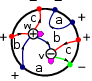} \end{ex}

Suppose that one takes a disk diagram and excises a neighborhood of each vertex. What remains is a colored, oriented 1-manifold embedded in the punctured disk.

Suppose that $f \co D^2 \to X^2$ is a map. It is \emph{represented} by a given disk diagram if \begin{itemize} \item For each $r \in \RC$, $f^{-1}(\hat{r})$ is the collection of vertices
labeled $r$. \item Each vertex labeled $r$ has a neighborhood which maps homeomorphically to the 2-cell $r$. The sign on the vertex is $+$ if the homeomorphism preserves orientation, and
$-$ if it reverses it. \item Let $Y$ denote the disk with those neighborhoods excised. Then $Y$ maps to $X^1$. The remainder of these criteria address the restricted map $f_Y \co Y \to
X^1$. \item For each $s \in \SC$, $f_Y^{-1}(\hat{s})$ is the 1-manifold colored $s$. \item Each connected 1-manifold colored $s$ in $Y$ has a tubular neighborhood mapping by projection to
a neighborhood of $\hat{s}$. The orientation of the manifold obeys the obvious rule. \end{itemize} Once again, these conditions imply that the remainder of the disk (i.e. $Y$ minus these
tubular neighborhoods) is sent to $X^1 \setminus \hat{\SC}$, and maps to the connected component corresponding to the label on each region.

\begin{defn} A map $D^2 \to X^2$ is \emph{pictorial} if it is represented by some disk diagram. A map $(S^2,\pt) \to X^2$ is \emph{pictorial} if the corresponding map $D^2 \to X^2$ is
pictorial, given by collapsing the boundary to $\pt$. This implies that the corresponding disk diagram is \emph{closed}, i.e. its boundary has the empty word. \end{defn}

Any disk diagram is clearly realized by some map $D^2 \to X^2$. Any two pictorial maps $D^2 \to X^2$ which agree on the boundary and have isotopic disk diagrams are clearly homotopic
relative to the boundary. In the definition of an enriched 3-complex, the attaching map of a 3-cell $z$ is assumed to be pictorial, and thus has a closed disk diagram. One can also define
the flip operation on disk diagrams, which inverts all the orientations.

Every map $D^1 \to X^1$ or $D^2 \to X^2$ is homotopic to a pictorial map (and if the boundary is already nice enough, this homotopy can be performed relative to the boundary). Any
3-complex is homotopy equivalent to a 3-complex with pictorial attaching maps. The choice of additional data needed to enrich a 3-complex is unique up to homotopy. Therefore, when studying
arbitrary maps from $D^2$ to arbitrary 3-complexes up to homotopy, it is sufficient to study pictorial maps from $D^2$ to enriched 3-complexes.

For more details, see Fenn \cite{Fenn}.

Henceforth, we will use the term \emph{Igusa diagram} to refer to any diagram (on the line, circle, or disk) constructed above. We only consider Igusa diagrams up to isotopy. We also use the
word \emph{symbol} to refer to a vertex in a disk diagram.

\subsection{Homotopy relations on diagrams}
\label{subsec-relations}

Consider a disk diagram with a sub-disk-diagram containing no symbols. This subdiagram represents a map $D^2 \to X^1$. There are two local transformations of diagrams which result in
homotopic maps $(D^2,\pa D^2) \to X^1$. These are called \emph{bridging} and \emph{removing circles}, and they can be applied to any $s \in \SC$ (we have omitted the labeling of regions).
We write the moves by placing an equal sign between the two diagrams. The transformation \eqref{removecircle} can also be applied with the other orientation.

\begin{equation} \ig{2}{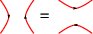} \label{bridge} \end{equation}
\begin{equation} \ig{2}{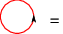} \label{removecircle} \end{equation}

\begin{claim} Any two symbol-less disk diagrams which yield homotopic maps $(D^2,\pa D^2) \to X^1$ are related by a sequence of \eqref{bridge} and \eqref{removecircle}. \end{claim}

If we allow symbols, there is a new local transformation of diagrams which results in homotopic maps $(D^2,\pa D^2) \to X^2$. It is called \emph{canceling pairs}, and can be applied to
any $r \in \RC$. In this relation, the orientations must be opposite and the tags must lie in the same region.

\begin{equation} \ig{2}{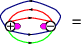} \label{cancelpair} \end{equation}

Together, \eqref{cancelpair}, \eqref{removecircle} and \eqref{bridge} are called the \emph{standard relations}.

\begin{exercise} Use \eqref{removecircle} and \eqref{bridge} to prove that \eqref{cancelpair} is equivalent to the local move

\begin{equation} \label{cancelpair2} \ig{2}{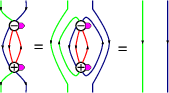} \end{equation} \end{exercise}

\begin{claim} Two diagrams with the same boundary represent relatively homotopic maps $(D^2,\pa D^2) \to X^2$ if and only if they are related by the standard relations. \end{claim}

The above claims are proven in \cite[\S 2.4]{Fenn}.

To construct a 3-complex from a 2-complex, one glues in a set $\ZC$ of oriented balls $D^3$ along maps $\pa D^3 \cong S^2 \to X^2$. The effect of adding a 3-cell $z \in \ZC$ to $X^2$ is
that it makes the corresponding closed diagram $(\pa z,\pt_z)$ nulhomotopic. The corresponding local move on disk diagrams would be to replace the diagram $(\pa z,\pt_z)$ with the empty
diagram, or vice versa. Note that a disk diagram always represents a map whose image lies in $X^2$, but this local move corresponds to a homotopy which passes through $X^3$. We typically do not bother to draw the tag corresponding to $\pt_z$ on such a disk diagram, because its location on the empty boundary is irrelevant.

Alternatively, one can also consider $\pa z \cong S^2$ as a union of two copies of $D^2$ along a common boundary $S^1$ (containing the marked point $\pt_z$). The two hemispheres would
represent two (possibly non-closed) diagrams with the same boundary, and the corresponding local move would be to replace one diagram with the other. Given a closed diagram, one can obtain
the two hemispheres by slicing the disk in half to form two disks, and taking the flip of one. The relation which replaces $\pa_z$ with an empty diagram and the relation which replaces one
hemisphere with another are equivalent modulo the standard relations.

\begin{ex} The following two relations, which could arise from the gluing of a 3-cell, are equivalent. The attaching maps $w$ and $v$ come from the previous examples. \igc{2}{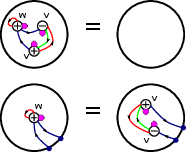}
\end{ex}

\begin{thm} Two diagrams with the same boundary represent relatively homotopic maps $(D^2,\pa D^2) \to X^3$ if and only if they are related by the standard relations and the new relations
imposed by $\ZC$. \end{thm}

\begin{remark} \label{remark:redundant3cell} If one glues in a new 3-cell along an attaching map $(\pa z,\pt_z)$ which is already nulhomotopic, then the new relation is clearly redundant.
In other words, any homotopy of maps $D^2 \to X^3$ which passes through $z$ could have instead avoided $z$ (though there may be no homotopy of homotopies). If two diagrams are homotopic,
it is easy to deduce from the standard relations that their flips are also homotopic. Therefore, after gluing in $z$, gluing in a 3-cell $\overline{z}$ along the flipped attaching map
will not affect the diagrammatic calculus. \end{remark}

\subsection{$2$-categorical language}
\label{subsec-nowin2cats}

We can also draw Igusa diagrams in the planar strip $\R \times [0,1]$ rather than the planar disk, and they will be called \emph{strip diagrams}. They represent (pictorial) maps
$(D^2,\pt,\pt) \to X^3$ with two marked points on the boundary. The same local moves as above will describe homotopy classes of such diagrams.

\begin{defn} Let $X^3$ be an (enriched) 3-complex, with 0-cells $\OC$, 1-cells $\SC$, 2-cells $\RC$, and 3-cells $\ZC$. We define a $2$-category $\Pi(X^3)_{\le 2}$ as follows. The objects
will be $\OC$. The $1$-morphisms will be generated by $s \co o_1 \to o_2$ and $s^{-1} \co o_2 \to o_1$, where $o_1$ (resp. $o_2$) is the 0-cell at the source (resp. target) of the oriented
edge $s$. Thus an arbitrary $1$-morphism is a compatible word in $\SC \cup \SC^{-1}$. The 2-morphisms $w_1 \to w_2$ between compatible words will be the set of strip diagrams constructed
with the symbols $r,r^{-1}$ for $r \in \RC$, modulo isotopy, the standard relations, and a relation for each $z \in \ZC$. Composition of 2-morphisms is given by vertical concatenation.
\end{defn}

\begin{remark} One can also phrase this definition in terms of generators and relations. The 2-morphisms are generated by oriented cups and caps for each $s \in \SC$, and by symbols $r$
and $r^{-1}$ for each $r \in \RC$. In addition to the standard relations and $\ZC$, one imposes certain ``isotopy relations." See Lauda \cite{LauSL2} for more details. \end{remark}

Note that oriented cups and caps give 2-morphisms $ss^{-1} \to \1$, etc. Relations \eqref{bridge} and \eqref{removecircle} prove that cups and caps form inverse isomorphisms. Similarly,
the symbol $r$ gives a map from $w(r) \to \1$, and $r^{-1}$ gives a map $\1 \to w(r)$. Relations \eqref{cancelpair2} and \eqref{cancelpair} prove that these are inverse isomorphisms.

This combinatorially-defined $2$-category encodes everything one needs to know about $\pi_{\le 2}(X^3)$. In particular, the previous results immediately imply this corollary.

\begin{cor} There is an obvious $2$-functor $\Pi(X^3)_{\le 2} \to \pi(X^3)_{\le 2}$, sending each object $o \in \OC$ to the corresponding point in $X^3$. This is a $2$-categorical
equivalence. \end{cor}

\subsection{Group presentations}
\label{subsec-presentations}

Let $\PC=(\SC,\RC)$ be a 2-presentation of a group $G$. The corresponding 2-complex is the \emph{Cayley complex} $X_{\PC}$, and is constructed in the familiar way. It has a single 0-cell,
a 1-cell for each $s \in \SC$, and a 2-cell for each $r \in \RC$, glued in the obvious fashion along its corresponding word (see also \cite[\S 1.2]{Fenn}). Note that $G \cong
\pi_1(X_{\PC})$, although the higher homotopy groups depend on the presentation chosen. Recall that a \emph{3-presentation} $\PC = (\SC,\RC,\ZC)$ is a 2-presentation of $G$ with a
collection $\ZC$ of 3-cells; we also denote the corresponding 3-complex $X_{\PC}$. We call the elements of $\ZC$ \emph{2-relations}.

\begin{ex} Suppose $G = \{e\}$ and $\PC = (\emptyset,\RC)$, where each element of $\RC$ is the empty word. Then $X_G$ will be a rosette of 2-spheres, one for each element of $\RC$. When
$\RC$ is a singleton so that $X_G \cong S^2$, the commutative group $\pi_2$ is isomorphic to $\Z$, based on a signed count of appearances of the relation symbol. \igc{1.5}{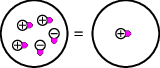}
\end{ex}

\begin{exercise} Suppose $G = \{e\}$ and $\PC= (\{a\},\{a\})$. Then $X_G \cong D^2$. Show explicitly that any two disk diagrams with the same boundary are equivalent. \end{exercise}

We can also construct a $G$-fold (universal) cover of $X_{\PC}$, which we will denote $\tX_\PC$. The 0-cells will be $\tOC = G$, the 1-cells will be $\tSC = \SC \times G$, the 2-cells will
be $\tRC = \RC \times G$, and so forth. Each 1-cell $(s,g)$ will go from $g$ to $gs$; by convention, edges correspond to right multiplication. Each 2-cell $(r,g)$ will be attached along
the edges corresponding to the word of $r$, beginning at the base point $g \in \tOC$. The 3-cell $(z,g)$ is glued into the closed diagram corresponding to $z$, with the outer region
labelled $g$, and the other regions labeled in the only consistent way. Clearly $\tX_\PC$ comes equipped with a free action of $G$ by left multiplication on cell names, and $\tX_\PC / G
\cong X_\PC$.

\begin{claim} $\pi_1(\tX_G) \cong \1$. Moreover, $\pi_n(\tX_G) \cong \pi_n(X_G)$ for all $n \ge 2$. \end{claim}

\begin{proof} This is immediate from the long exact sequence associated to the covering map. \end{proof}

\begin{remark} It is somewhat presumptuous to assume that, given $\PC$, one knows what $G$ is, or even how big $G$ is. While $X_\PC$ can be constructed explicitly without the set of
elements of $G$, $\tX_\PC$ can not be. When $G$ is infinite so is $\tX_\PC$, but it is locally finite so that usual topological intuition applies. \end{remark}

We use the following conventions for Igusa diagrams of group presentations. We do not bother to label the regions of a diagram for $X_\PC$ with the unique element of $\OC$. Given a disk
diagram for $\tX_\PC$, the label $g \in G = \tOC$ of a single region will determine the label of every other region. Moreover, having fixed a label on a single region, the color $(s,g) \in
\tSC$ of any edge is determined only by the color $s \in \SC$, and similarly for symbols $(r,g) \in \tRC$. We omit the redundant data, coloring edges only by $s \in \SC$ and naming symbols
by $r \in \RC$. Thus a disk diagram for $\tX_\PC$ is the same data as a disk diagram for $X_\PC$ with an arbitrary choice of label $g \in G$ for a single chosen region. Postcomposing a map $D^2 \to
\tX_\PC$ with the quotient map $\tX_\PC \to X_\PC$ corresponds to forgetting the label on that region.

The following example will be crucial in chapter \ref{sec-modifiedFenn}.

\begin{ex} \label{involutionexample} Let $\PC = (\{s\},\{s^2\},\{z\})$ be a presentation for the group $G = \Z/2\Z$. The 3-cell $z$ is glued in along the picture:
\igc{1.75}{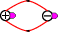} Then the 2-skeleton of $X_\PC$ will be $\RM \PM^2$, and $X_\PC$ will be $\RM \PM^3$. In particular, $\pi_2(X_\PC)$ is trivial.

To construct $\tX_\PC$, we take two points, add two edges to get $S^1$, add two disks to get $S^2$, and add two 3-cells to get $S^3$. This is the 3-skeleton of $S^\infty \cong EG$ in its
usual construction.

Let $Y$ denote the 2-skeleton of $\tX_\PC$ with only a single 3-cell added, so that $Y \cong D^3$. We think of $Y$ as the ``2.5-skeleton" of $EG$. To obtain $\tX_\PC$ from $Y$, one
attaches a 3-cell which is redundant in the sense of Remark \ref{remark:redundant3cell}. Therefore, $Y$ and $\tX_{\PC}$ have the same category $\Pi_{\le 2}$. Of course, $Y$ does not admit
a free action of $G$, but it has other advantages. For instance, $Y$ deformation retracts to a pole between the two 0-cells, which is the (completed) dual Coxeter complex of $G$. \end{ex}

\section{Modified Fenn Diagrams}
\label{sec-modifiedFenn}

The Coxeter presentation has a number of natural symmetries, and we wish to exploit them in order to simplify our diagrammatic description of $\Omega W$. In this chapter we develop some
general machinery which yields simpler diagrammatics for special kinds of 3-presentations.

Suppose that $(\SC,\RC)$ is a group presentation, where $s^2 \in \RC$ for some $s \in \SC$. There is a particular 3-cell one can glue in, which will cause $s$ and $s^{-1}$ to be
canonically isomorphic, and this allows us to ignore the orientations on the strands colored $s$ in diagrams for $\tX_\PC$. Heuristically, these $s$-unoriented diagrams depict $\Pi_{\le
2}$ for some deformation retract of a ``2.5-skeleton" of $\tX_{\PC}$, as in Example \ref{involutionexample}. In similar fashion, we describe a modification adapted to rotational and flip
symmetries in relations, such as in the braid relation.


\subsection{Modified diagrammatics for involutions}
\label{subsec-involutions}

First let us consider diagrams for $\PC = (\{s\},\{s^2\})$, so that $X_\PC \cong \RP^2$ and $\tX_\PC \cong S^2$. The relation $s^2$ allows us to draw bivalent vertices which look like
this. \igc{1.5}{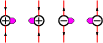} The sign on the symbol is determined by the orientations of the strands, but the location of the tag is not. Therefore, the bivalent vertex gives two natural maps
$s \to s^{-1}$, depending on the placement of the tag, and two natural maps the other direction.

Using \eqref{cancelpair} and \eqref{cancelpair2} we have

\begin{equation} \ig{1.5}{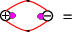} \label{RP2cancel} \end{equation}

\begin{equation} \ig{1.5}{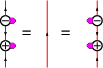} \label{RP2isom1} \end{equation}

We now add a $3$-cell $z_s$ to obtain the higher presentation $\PC=(\{s\},\{s^2\},\{z_s\})$, and temporarily write $\tX=\tX_\PC$ and $X=X_\PC$. The new $3$-cell is meant to
kill $\pi_2(\RP^2)$, and in $\tX$ to kill $\pi_2(S^2)$. If we attach two bivalent vertices together so that the tags do not cancel, this represents the map that $z_s$ is glued
into. Thus we have a new relation:

\begin{equation} \ig{1.5}{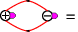} \label{RP2notcancel} \end{equation}

Splitting a new $3$-cell into hemispheres, we obtain the equivalent relations:

\begin{equation} \ig{1.75}{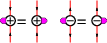} \label{RP2tagswitch} \end{equation}

We may introduce a new symbol: a bivalent vertex without a tag. This symbol is set equal to the bivalent vertex with either placement of the tag. Now \eqref{RP2isom1} becomes

\begin{equation} \ig{1.5}{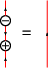} \label{RP2isom2} \end{equation} 

Thus the bivalent vertex gives a (canonical) isomorphism between $s$ and $s^{-1}$.

\begin{exercise} Any closed diagram for $X_\PC$ is equal to the empty diagram. In other words, $\pi_2(X_\PC)=1$. The topological statement is obvious: find a diagrammatic proof.
\label{noclosedforRP2now} \end{exercise}

Now let $\PC = (\SC,\RC,\ZC)$ be any 3-presentation, and let $s \in \SC$ be such that $s^2 \in \RC$ and $z_s \in \ZC$ for $z_s$ as above. The same arguments as above show that we may
ignore the tag on the bivalent vertex associated to $s$. Moreover, the bivalent vertices of either sign form inverse isomorphisms between $s$ and $s^{-1}$, and we wish to use them to
canonically identify the two objects. Given any Fenn diagram we may forget the orientation data associated to $s$ to get an \emph{$s$-unoriented Fenn diagram}. In particular, each relation
has an $s$-unoriented symbol. Let $\overline{s}$ denote a point $s$ without an orientation.

\begin{defn} Suppose that $\PC=(\SC,\RC,\ZC)$ is a 3-presentation containing $(s,s^2,z)$. Let $\Pi(X_\PC,s)_{\le 2}$ be the 2-category with a single object, defined as follows. The
1-morphisms are generated by $\SC' \cup (\SC')^{-1} \cup \{\overline{s}\}$, where $\SC' = \SC \setminus \{s\}$. The 2-morphisms are generated by the $s$-unoriented symbols of $\RC' \cup
(\RC')^{-1}$, for $\RC' = \RC \setminus \{s^2\}$, and thus correspond to $s$-unoriented disk diagrams. The relations are generated by $\ZC \setminus \{z_s\}$, as well as the usual Fenn
relations for oriented parts of the diagram and the \emph{unoriented Fenn relations} for $s$:

\begin{equation} \ig{1.5}{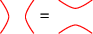} \label{bridgeunor} \end{equation}
\begin{equation} \ig{1.5}{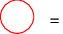} \label{removecircleunor} \end{equation}
\begin{equation} \ig{1.75}{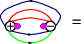} \label{cancelpairunor} \end{equation}

\end{defn}

There is a natural 2-functor $\Pi(X_\PC)_{\le 2} \to \Pi(X_\PC,s)_{\le 2}$. It sends both $s$ and $s^{-1}$ to $\overline{s}$. It sends the bivalent vertex corresponding to $s^2$ to the
identity map of $\overline{s}$. To every disk diagram without bivalent vertices, it forgets the orientation data associated to $s$. It is easy to show that this 2-functor is an
equivalence. Given any $s$-unoriented disk diagram, and any choice of orientations of $s$ on the boundary, one may choose a disk diagram by placing orientations on $\overline{s}$-strands
willy-nilly, and adding bivalent vertices whenever necessary for consistency. While there are multiple such diagrams, they are all equal in $\Pi(X_\PC)_{\le 2}$.

If $\PC$ contains $(s,s^2,z_s)$ for multiple distinct involutions in $\SC$, there is no obstruction to forgetting the orientations on multiple colors at once. We write
$\Pi(X_\PC)^{\unor}_{\le 2}$ for the 2-category where every such orientation is ignored.

The case of $\tX_\PC$ can be treated in the same way. One must glue in a copy of $z_s$ for every possible region labeling. As before, diagrams for $\tX_\PC$ will be $s$-unoriented disk
diagrams with a label in a single region.

\begin{remark} \label{defretract1} Here is a heuristic topological understanding of unoriented diagrams, at least for $\tX_\PC$. As in Example \ref{involutionexample}, let $Y\cong D^3$ be
the 2.5-skeleton of $S^\infty$, which has two 0-cells $1$ and $s$, two 1-cells $(s,1)$ and $(s,s)$, two 2-cells $(s^2,1)$ and $(s^2,s)$, and a single 3-cell $(z_s,1)$. One can construct a
new ``unoriented" cell complex $Y^{\unor}$, consisting of two 0-cells $1$ and $s$, and a single unoriented 1-cell $\overline{s}$ between them. We think of $\overline{s}$ as a pole inside
$Y \cong D^3$. Clearly $Y^{\unor} \subset Y$ is a deformation retract, under a retract sending both edges $(s,1)$ and $(s,s)$ to $\overline{s}$.

Similarly, suppose that $(s,s^2,z_s) \subset \PC$ for a general 3-presentation. After constructing $\tX_\PC^1$, one can repeat the above construction for each coset $\{x,xs\} \in G$ to
obtain a 3-complex $Y_1$ which deformation retracts to a 1-complex $Y_1^{\unor}$. In $Y_1^{\unor}$, $x$ and $xs$ are connected by a single edge $(\overline{s},x)$. The attaching maps of
other 2-cells in $\RC$ can be deformed to lie on $Y_1^{\unor}$, and similarly for the other 3-cells in $\ZC$, yielding a deformation retract $Y^{\unor}$ of a 2.5-skeleton $Y$ of $\tX_\PC$.
We think of unoriented diagrams as describing maps to $Y^{\unor}$ (even though Fenn diagrams for $Y^{\unor}$ are actually quite different). There is no reasonable $\Z/2\Z$ action on $Y$ or
$Y^{\unor}$ whose quotient has $\pi_1 = G$, so we do not use this heuristic when thinking about $X_\PC$, only $\tX_\PC$. \end{remark}

\subsection{Rotational invariance and flip invariance}
\label{subsec-rotations}

When a relation does not have rotational invariance, there is no need to keep track of the tag on the corresponding symbol in a Fenn diagram. The location of the tag can be deduced from
the edge coloring. When a relation does have rotational invariance, the tag is not redundant. However, if an appropriate 3-cell is glued in, all possible locations of the tag will be set
equal, and the tag will become redundant. An exactly analogous procedure will work to make the sign on a symbol redundant when a relation has flip invariance. We will use specific examples
to illustrate general principles, because it is hard to draw a general example.

First consider $\PC = (\{r,g,b\},\{w=rgbrgbrgb\})$. The symbol for $w$ can be rotated by 120 degrees and 240 degrees to give a morphism with the same boundary. This is a different morphism
because the tag is in the wrong place.

\igc{1.5}{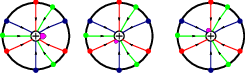}

If we set two of these to be equal by gluing in a 3-cell, then the third will be equal as well. In general, if $w$ is invariant under rotation by $\theta$ then setting $w$ equal to
$\theta(w)$ will also set it equal to $n\theta(w)$ for any $n \in \Z$. The 3-cell $z_w$ would be glued in along the following closed diagram, which is a ``mismatched pair."

\igc{1.5}{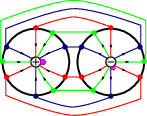}

Once this 3-cell is glued in, one need not draw the tag on this symbol any longer. There are only three valid locations for the tag (it must be before $r$
and after $b$), and they all give equal 2-morphisms.

In the previous section, we had to construct a new 2-category which equated two canonically isomorphic objects. In this section, we are not changing the category, but are merely using a
notational convenience, using one symbol to represent several distinct symbols which happen to be equal.

The case of $\tX_\PC$ can be treated in the same way. One must glue in a 3-cell as above for every possible region labeling.

\begin{remark} \label{defretract2} As in Remark \ref{defretract1}, there is a topological heuristic for the new diagrammatic calculus. Suppose that $\PC$ contains
$(\{r,g,b\},\{w=rgbrgbrgb\}, z_w)$ as above. For any $x \in G$ there are three different $2$-cells being glued to the same $S^1 \subset \tX_\PC$: $(w,x)$, $(w,xrgb)$ and $(w,xrgbrgb)$.
Fixing the same base point in $S^1$ for all three, they are $(w,x)$, $(\theta(w),x)$ and $(\theta^2(w),x)$. One can visualize this part of the 2-skeleton as a stack of pancakes, glued
together along their rim. The 3-cell $(z_w,x)$ fills in the gap between the first two pancakes, while the 3-cell $(z_w,xrgb)$ fills in the gap between the second and third pancakes. With
these two 3-cells glued in, the result is a copy of $D^3$. Therefore the last 3-cell $(z_w,xrgbrgb)$ would be redundant, and we need not glue it in. Ignoring this 3-cell (for each $x$) one
obtains the ``2.66-skeleton" of $\tX_\PC$ (we continue to call it the 2.5-skeleton), and it deformation retracts to a central $D^2$ pancake-shaped slice. This central slice is what the
tagless symbol is meant to represent.

The reader can deduce the rest of the analogy. Unlike Remark \ref{defretract1}, replacing $\tX_\PC$ with the deformation retract of its 2.5-skeleton does not change the 1-skeleton, which
is why one need not change the objects in the category. \end{remark}

Now consider $G=(\{r,g\},\{w=rgg^{-1}r^{-1}\})$. It lacks any rotational symmetry, but it does have a flip symmetry: the symbol for $w$ and some rotation of the opposite orientation of $w$
have the same boundary. We can glue in a 3-cell $z_w$ along a ``mismatched pair."

\igc{1.5}{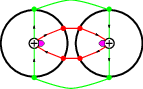}

Once this 3-cell is glued in, one need not keep track of the sign on the symbol any longer.

Of course, the relation $w \in \pi_1(X_\PC)$ is already nulhomotopic even in $X_\PC^1$, as any relation with flip symmetry will be! This restricts the notion of flip symmetry to
unusual presentations.

Flip symmetry becomes more interesting for unoriented Fenn diagrams. Suppose that $\PC = (\{r,g,b\},\{r^2,g^2,b^2,rgbrbg=w\},\ZC)$ and that $\ZC$ contains the 3-cells which allow for
unoriented diagrams as in the previous section. The unoriented symbol for $w$ has no rotational symmetry, but it does have flip symmetry.

\igc{1.5}{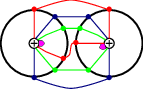}

Without this 3-cell added as a relation, a diagram with boundary $rgbrbg$ could be either $w_+$ or a rotated $w_-$. This 3-cell would set them equal. For this example it is not terribly
meaningful to say that we can remove the $\pm$ decoration on $w$, because one must keep the tag, and the sign can be deduced from the tag. The following example combines all three modifications, and gives a situation where removing the sign does have a noticeable effect.

\begin{ex} Consider the Coxeter presentation $(\{s,t\},\{s^2,t^2,stst^{-1}s^{-1}t^{-1}\})$. Now, glue in the 3-cells for each generating involution, so we may work with unoriented
diagrams. Then glue in a 3-cell for rotational invariance. At this point, the following two diagrams do not represent the same 2-morphism (tags included for clarity).

\igc{1.5}{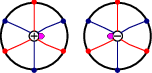}

Gluing in one more 3-cell for flip symmetry, we can ignore the sign on the symbol, and draw the 2-morphism unambiguously as a 6-valent vertex.

\igc{1.5}{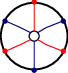} \end{ex}

\begin{ex} \label{pancakeexample} Let $\PC$ be the 3-presentation of the previous example. What is $\tX_\PC$? It is more complicated than it looks, because each cell appears 6 times, once
for each element of $W$.

We begin with six 0-cells. Instead of labeling them by elements of $W$, let us mentally arrange them as a hexagon and label them by numbers modulo 6, as on a 6-hour clock. We then glue in
twelve 1-cells (six for $s$ and six for $t$). These connect each neighboring pair of 0-cells with two edges, yielding six copies of $S^1$ welded into a hexagonal loop. Next, the relations
in $\QC$ attach twelve 2-cells (six of each), two glued into each copy of $S^1$. This yields six copies of $S^2$ welded into a loop. Then we glue in the ``orientation 3-cells", twelve of
them, two glued into each copy of $S^2$. This yields six copies of $S^3$ welded into a loop. However, six of the twelve 3-cells are redundant: after a 3-cell has turned $S^2$ into $D^3$,
the other 3-cell will have a nulhomotopic attaching map. Ignoring these six redundant 3-cells, we have six copies of $D^3$ welded into a loop. This space deformation retracts to the
1-skeleton of $\tCox$ (to be defined in \S \ref{subsec-dualcoxetercomplex}), which is just (the boundary of) a hexagon.

Now we glue in six 2-cells, for the braid relation. It may help the reader to think that three of these are glued in clockwise, and three counter-clockwise, but the underlying topological
space does not care about such intricacies. The result is a stack of 6 hexagonal pancakes (labeled by the numbers modulo 6), glued along their boundary.

Next we glue in the six 3-cells for rotational invariance. One such 3-cell forms a cobordism between pancake 0 and pancake 2, another between 2 and 4, and another between 4 and 0; the
three remaining 3-cells go between pancakes 5 and 3, pancakes 3 and 1, and pancakes 1 and 5. Clearly one of the even 3-cells and one of the odd 3-cells is redundant.

Finally, we glue in six 3-cells for flip invariance. One such 3-cell forms a cobordism between pancake 0 and pancake 3, another between 1 and 4, another between 2 and 5, another again
between 3 and 0, and so forth. At this point, after gluing in one flip 3-cell, the remaining ones are redundant.

Ignoring the redundant 3-cells, we glue in five 3-cells to fill the gaps between the six pancakes. The result is a big blob of pancake batter, which clearly flattens into a single solid
hexagon. In other words, this space (minus the redundant 3-cells) deformation retracts to the 2-skeleton of $\tCox$, a solid hexagon. \end{ex}

\part{Coxeter groups and braid groups}
\label{pt:coxeter}

\section{Coxeter groups and topology}
\label{sec-coxbackground}

In this chapter we give some background information on Coxeter groups, their Artin braid groups, and some associated topological spaces.

\subsection{Coxeter groups}
\label{subsec-coxetergroups}

Fix a set $\SC$, and for each pair $s \ne t \in \SC$ fix an element $m_{st} \in \Z_{\ge 2} \cup \{\infty\}$. The Coxeter group $W$ is defined by its \emph{Coxeter presentation} $(\SC,\QC
\cup \BC)$, where the quadratic relations are \[\QC = \{s^2\}_{s \in \SC}\] and the braid relations are \[\BC = \{b_{s,t}\}_{s \ne t \in \SC} \text{ for } b_{s,t} =
\ubr{sts\ldots}{m_{st}}\ubr{\ldots t^{-1}s^{-1}t^{-1}}{m_{st}}.\] There is no braid relation when $m_{st}=\infty$. There is only one braid relation for each pair $s,t \in \SC$; we will not
redundantly use both $b_{s,t}$ and $b_{t,s}$. The corresponding Artin braid group $B_W$ has presentation $(\SC,\BC)$. We let $B_W^+ \subset B_W$ denote the monoid of \emph{positive
braids}, which is the monoid with the same presentation $(\SC,\BC)$.

We assume that $\SC$ is finite, though this is not strictly necessary for our arguments. We let $r = |\SC|$ be the \emph{rank} of $W$. Let $\ell$ denote the length function.

For a subset $I \subset \SC$, there is a parabolic subgroup $W_I \subset W$ generated by $s \in I$. It is also a Coxeter group, with presentation $(I,\QC_I \cup \BC_I)$. When there exists
a partition $\SC = I_1 \coprod I_2$ such that $m_{st}=2$ for all $s \in I_1$ and $t \in I_2$, then $W \cong W_{I_1} \times W_{I_2}$, and we say that $W$ is \emph{reducible}. When $W_I$ is
finite, we say that $I$ is \emph{finitary}, and we let $w_I$ denote the longest element of $W_I$.

A Coxeter group of rank 2 is determined by $m=m_{st}$, and is said to be of type $I_2(m)$. It is finite unless $m = \infty$. The group $I_2(2)$ is the reducible group $A_1 \times A_1$.
There is a classification of all finite Coxeter groups. The finite Coxeter groups of rank 3 are types $A_3$, $B_3$, $H_3$ and the reducible types $A_1 \times I_2(m)$ for $m < \infty$.

For an element $w \in W$, we will use an underline $\ul{w} = s_1 s_2 \cdots s_d$ to indicate an expression for $w$ in terms of $\SC$. If we need to differentiate between two expressions
for $w$ we will write $\ul{w}_1$ and $\ul{w}_2$. We say that $\ul{w}$ is \emph{reduced} if $d=\ell(w)$. Given $w \in W$, a choice of reduced expression $\ul{w}$ will also yield an element
$\tilde{w}$ of $B_W^+$, independent of the reduced expression chosen. We call this the \emph{positive lift} of $w$ to $B_W$.

See Humphreys \cite{Humphreys} for more details.

\subsection{The Coxeter complex}
\label{subsec-coxetercomplex}

To a Coxeter system $(W,\SC)$ one may associate a simplicial complex, the \emph{Coxeter complex} $|(W,\SC)|$, as follows: \begin{enumerate} \item Choose an arbitrary total order on $\SC$.
\item Color the $r$ faces of the $(r-1)$-simplex by $\SC$, matching the lexicographic order on faces to the total order on $\SC$; call the resulting simplex $\Delta$. \item Take one copy
$\Delta_w$ of $\Delta$ for each $w \in W$. \item Glue $\Delta_w$ to $\Delta_{ws}$ along the face colored by $s$, for all $w \in W$ and $s \in \SC$. There is only one possible gluing which
preserves the orientation. \end{enumerate} The result is a connected $(r-1)$-dimensional simplicial complex with simplices of maximal dimension labelled by $W$ and codimension one
simplices (or \emph{walls}) colored by $\SC$. Moreover, $W$ acts on $|(W,\SC)|$ by automorphisms preserving the coloring of walls.

If one chooses a different total order on $\SC$, one obtains the same complex with different orientations. We do not care about the simplicial orientations in the Coxeter complex (or in
any of the other complexes we construct in this chapter; the avid reader may fill in the details). In the examples below we will draw orientations on the walls; these orientations have
nothing to do with the simplicial orientations, but instead record the Bruhat order on $W$.

\begin{ex} Let $W$ be a finite subgroup of the orthogonal group $O(V)$ of a Euclidean vector space $V$ of dimension $r$, and assume $W$ is generated by reflections and acts irreducibly.
Let $T$ denote the subset of $W$ of elements which act as reflections on $W$. Consider the space \[ U := V - \bigcup_{t \in T} V^t \] obtained from $V$ by deleting all reflecting
hyperplanes. Then $W$ acts simply transitively on the connected components of $U$. If one fixes a connected component $C$ of $U$, then $\overline{C}$ is a simplicial cone, and $(W,\SC)$ is
a Coxeter system of rank $r$, where $\SC$ denotes the set of reflections in the walls of $\overline{C}$. If one intersects $\overline{C}$ with the unit sphere in $V$ then one obtains a
closed subset $\Delta$ homeomorphic to an $(r-1)$-simplex, whose faces are colored by $\SC$. The $W$-translates of $\{ w\Delta \; | \; w \in W \}$ give a triangulation of the unit sphere,
giving a realization of the Coxeter complex. In fact, all Coxeter complexes associated to finite Coxeter systems can be realized in this way. \end{ex}

\begin{ex} As for any finite rank 3 Coxeter system, the Coxeter complex for $A_1 \times A_1 \times A_1$ is a triangulation of the sphere. The triangles are labeled by $w \in W$. The
triangle closest to the reader is labeled with the identity, and the triangle furthest is the longest element. We place orientations on edges such that going from the left side of an edge
to the right side will increase the length of $w \in W$ by $1$.

 \begin{equation*} \begin{tikzpicture}[scale=1.5] \def\ct{0.6cm}

	\coordinate (d1) at (0:1cm);
	\coordinate (d2) at (60:1cm);
	\coordinate (d3) at (120:1cm);
	\coordinate (d4) at (180:1cm);
	\coordinate (d5) at (-120:1cm);
	\coordinate (d6) at (-60:1cm);

	\coordinate (e1) at (30:\ct);
	\coordinate (e2) at (90:\ct);
	\coordinate (e3) at (150:\ct);
	\coordinate (e4) at (-150:\ct);
	\coordinate (e5) at (-90:\ct);
	\coordinate (e6) at (-30:\ct);

	\node[circle,draw,inner sep = 0mm, minimum size=2mm] (c1) at (30:1cm) {};
	\node[circle,draw,inner sep = 0mm, minimum size=2mm] (c2) at (90:1cm) {};
	\node[circle,draw,inner sep = 0mm, minimum size=2mm] (c3) at (150:1cm) {};
	\node[circle,draw,inner sep = 0mm, minimum size=2mm] (c4) at (-150:1cm) {};
	\node[circle,draw,inner sep = 0mm, minimum size=2mm] (c5) at (-90:1cm) {};
	\node[circle,draw,inner sep = 0mm, minimum size=2mm] (c6) at (-30:1cm) {};

	\draw[color=red] (c5) .. controls (e6) .. node[sloped]{$>$} (c1);
	\draw[color=red] (c1) .. controls (d2) .. node[sloped]{$<$} (c2);
	\draw[color=red!25] (c2) .. controls (e3) .. node[sloped]{$<$} (c4);
	\draw[color=red] (c4) .. controls (d5) .. node[sloped]{$>$} (c5);

	\draw[color=blue] (c1) .. controls (e2) .. node[sloped]{$<$} (c3);
	\draw[color=blue] (c3) .. controls (d4) ..  node[sloped]{$>$} (c4);
	\draw[color=blue!25] (c4) .. controls (e5) .. node[sloped]{$>$} (c6);
	\draw[color=blue] (c6) .. controls (d1) .. node[sloped]{$>$} (c1);

	\draw[color=black] (c3) .. controls (e4) .. node[sloped]{$>$} (c5);
	\draw[color=black] (c5) .. controls (d6) .. node[sloped]{$>$} (c6);
	\draw[color=black!25] (c6) .. controls (e1) .. node[sloped]{$<$} (c2);
	\draw[color=black] (c2) .. controls (d3) .. node[sloped]{$<$} (c3);

	\end{tikzpicture}
	\end{equation*}
\end{ex}

Given $x, y \in W$, a \emph{gallery from $x$ to $y$} in the Coxeter complex $|(W,S)|$ is a path between the simplices corresponding to $x$ and $y$, which does not meet any simplex of
codimension $\ge 2$. We regard two galleries as equivalent if they visit the same simplices in the same order. A gallery from $x$ to $y$ is \emph{minimal} if it crosses the least number of
walls amongst all galleries from $x$ to $y$. Giving a gallery from $x$ to $y$ is the same thing as giving an expression for $x^{-1}y$. Indeed, a gallery is determined uniquely by the
ordered list of walls crossed in the path from $x$ to $y$. A gallery from $x$ to $y$ corresponding to an expression $st\cdots u$ for $x^{-1}y$ is minimal if and only if $st \cdots u$ is
reduced.

\subsection{The Dual Coxeter complex}
\label{subsec-dualcoxetercomplex}

For our purposes it will be more convenient to use the \emph{dual Coxeter complex}, which is the CW-complex $|(W,\SC)|^\vee$ dual to $|(W,\SC)|$. It has a 0-cell for each $w \in W$,
and a gallery in the Coxeter complex corresponds to a path in the 1-skeleton of the dual Coxeter complex.

Let $C$ be any face of codimension $k<r$ in $|(W,\SC)|$. One can label $C$ by the rank $k$ subset $I \subset \SC$, consisting of the colors on the walls which contain that face. Then there
is a face labeled by $I$ if and only if $I$ is finitary. Moreover, the $(r-1)$-simplices containing such a face $C$ are labeled by elements of $W$ forming a coset in $W / W_I$.

Hence one can construct $|(W,\SC)|^\vee$ as follows: \begin{enumerate} \item Take a 0-cell for each $w \in W$. \item Attach a 1-cell from $x$ to $xs$, when $xs>x$. \item Attach
a 2-cell between the two minimal galleries from $x$ to $xw_{s,t}$, when $m_{s,t}$ is finite and $x$ is a minimal length coset representative. \item $\dots$ \end{enumerate} Let us
elaborate upon the inductive step. Fix any coset $C$ in $W/ W_I$ for $I$ finitary of rank $k$. Consider the cells whose closure only contains 0-cells corresponding to elements in
$C$. After the $k-1$-st step, the union of these cells will be homeomorphic $S^{k-1}$. The $k$-th step is to glue in a $k$-cell and obtain $D^k$ instead. As $|(W,\SC)|$ is $(r-1)$-dimensional, this process ends after $(r-1)$ steps.

One can also form the \emph{completed dual Coxeter complex}, which includes the $r$-th step above. We denote it by $\tCox$. It differs from $|(W,\SC)|^\vee$ in
a single $r$-cell when $W$ is finite, and does not differ otherwise. In the finite case, $\tCox$ gives a CW-complex structure for the unit ball in Euclidean space, rather than the unit
sphere.

\begin{ex} When $W$ is a finite dihedral group of size $2m$, $\tCox$ is the solid $2m$-gon. \end{ex}

\begin{exercise} Suppose that $W = W_1 \times W_2$ is a product of two other Coxeter groups. Show that $\tCox_W \cong \tCox_{W_1} \times \tCox_{W_2}$, compatibly with the CW structure.
\end{exercise}

\begin{prop} The completed dual Coxeter complex $\tCox$ is contractible. \end{prop}

\begin{proof} By the exercise above, we may assume that $W$ is irreducible. When $W$ is infinite, the result follows from the contractibility of the Coxeter complex (see e.g. \cite[Theorem 4.127]{Brown}). When $W$ is finite, the completed dual Coxeter complex is a unit ball. \end{proof}

The (completed) dual Coxeter complex $\tCox$ does have an action of $W$, which acts by left multiplication on 0-cells. However, this action is not free, and the quotient does not inherit a
nice CW-complex structure. The dihedral group acting on the regular $2m$-gon provides a familiar example. Thus $\tCox$ does not provide a good CW-complex model for $EG$. Instead, the
(3-dimensional) model we construct in chapter \ref{sec-coxdiag} will contain (the 3-skeleton of) $\tCox$ as a deformation retract.

%
%
%





\subsection{The Salvetti complex}
\label{subsec-Kpi1}

The completed dual Coxeter complex has one $k$-cell for each pair $(I,C)$, where $I$ is finitary of rank $k$, and $C$ is a coset of $W/W_I$. Suppose we place an equivalence relation on
$\tCox$, identifying any two $k$-cells $(I,C)$ and $(I,C')$. The quotient is still a CW-complex, having a single $k$-cell for each finitary $I \subset \SC$. For instance, there is a single
0-cell, a 1-cell for each $s \in \SC$, and a 2-cell for each pair $s \ne t \in \SC$ with $m_{s,t}<\infty$. We call this CW-complex $|B_W|$.

Similarly, one can construct a $W$-fold cover of this CW-complex, called the \emph{Salvetti complex} $\Sal$. It has one $k$-cell for each pair $(I,w)$, with $w \in W$ and $I \subset \SC$
finitary of rank $k$. The $k$-cell $(I,w)$ is glued in such a way that it contains 0-cells labeled by $wu$ for $u \in W_I$.

Note that the Salvetti complex is different from $\tCox$, despite having the same 0-cells, and $|B_W|$ is different from the quotient of $\tCox$ by the action of $W$ described above.

\begin{ex} Consider type $A_1$. Then $\tCox \cong D^1$ is an interval connecting two 0-cells $1$ and $s$. The complex $\Sal \cong S^1$ has two $1$-cells connecting the 0-cells $1$ and
$s$. The quotient $\tCox/W$ is also an interval, folded in half. Meanwhile $|B_W| \cong S^1$ identifies the endpoints of the interval, or wraps the Salvetti complex in half. \end{ex}

We have already discussed the $K(\pi,1)$-conjecture in some detail in the introduction \S\ref{subsec-Kpi1intro}. It is clear that $\pi_1(|B_W|) \cong B_W$. The $K(\pi,1)$-conjecture states
that all higher homotopy groups vanish; the $K(\pi,1)$-conjecturette states that $\pi_2(B_W)=0$. The $K(\pi,1)$-conjecturette is known for all Coxeter groups $W$, thanks to work of
Digne-Michel \cite{DigneMichel}.

For more information on the $K(\pi,1)$-conjecture and a list of cases where it is known, see the survey paper \cite{Paris2}.

\section{Diagrammatics for Braid groups}
\label{sec-braiddiag}

In \S\ref{subsec-Kpi1} we have already described the CW-complex $|B_W|$. We now seek to describe $\pi(|B_W|)_{\le 2}$ diagrammatically.

\begin{defn} Let $\Br^\diag$ denote $\Pi(X_\PC)$ for the 3-presentation $\PC=(\SC,\BC,\ZC)$ below. The presentation $(\SC,\BC)$ agrees with the presentation of the braid group given in
\S\ref{subsec-coxeterintro}. Therefore, an object in $\Pi(X_\PC)$ is a word in the letters $\SC \cup \SC^{-1}$. The morphisms are generated by oriented cups and caps, as well as
$2m_{st}$-valent vertices as pictured below, whenever $m_{st}<\infty$ for the two colors present.

\igc{1}{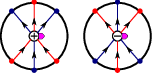}

These morphisms satisfy the Fenn relations:

\begin{equation} \ig{2}{bridge} \label{bridgeBraid} \end{equation}
\begin{equation} \ig{2}{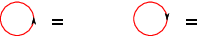} \label{removecircleBraid} \end{equation}
\begin{equation} \ig{2}{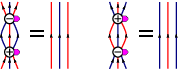} \label{6valentFenn} \end{equation}

\begin{remark} There is only one 2-cell for each pair $s \ne t \in \SC$ with $m_{s,t}<\infty$. The two different kinds of $2m$-valent vertices are the two orientations of the corresponding
symbol. Both the tag and the orientation on the symbol can be determined from the coloring and orientation on the strands, so we do not draw them in our diagrams henceforth. \end{remark}

In addition, for any three colors forming a finite parabolic subgroup, there is a single 3-cell in $\ZC$. The corresponding relation is the generalized Zamolodzhikov equation, given in
\eqref{Zamrelationsintro}. \end{defn}

By now, it is clear that Theorem \ref{mainthmbraid} is equivalent to the $K(\pi,1)$-conjecturette, and is thus proven.

\section{Diagrammatics for Coxeter groups}
\label{sec-coxdiag}

Let $(W,\SC)$ be a Coxeter group, with the usual presentation $(\SC,\QC \cup \BC)$. Let $\tCox$ be its completed dual Coxeter complex.

\begin{defn} A \emph{standard diagram} for $W$, will be a diagram with unlabeled regions, unoriented edges colored by $s \in \SC$, and (untagged, unoriented) $2m$-valent vertices which
alternate between edges colored $s$ and $t$ for which $m_{st}=m < \infty$. A \emph{labeled standard diagram} is a standard diagram with a single region labeled by an element of $W$.
\end{defn}

As noted previously, it is equivalent to give a label in $W$ for a single region, and to consistently label each region by an element of $W$, such that two regions separated by an edge $s$
differ by that element in $W$.

\begin{defn} Let $\WC^\diag$ denote the monoidal category whose objects are generated by $s \in \SC$, and whose morphisms are given by standard diagrams modulo isotopy and the following
relations (the Fenn relations and the Zamolodzhikov relations).
	
\begin{equation} \ig{1.5}{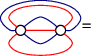} \end{equation}
\begin{equation} \ig{1.5}{removecircleunor} \end{equation}
\begin{equation} \ig{1.5}{bridgeunor} \end{equation}

\begin{equation} \ig{1}{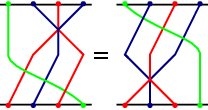} \end{equation}
\begin{equation} \ig{1}{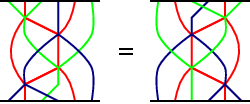} \end{equation}
\begin{equation} \ig{1}{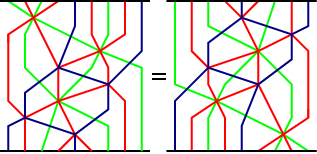} \end{equation}
\begin{equation} \ig{1}{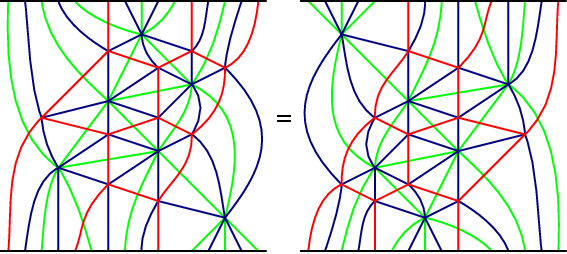} \end{equation}

Let $\tWC^\diag$ denote the 2-category whose objects are elements of $W$, and whose 2-morphisms are labeled standard diagrams, modulo the relations above. \end{defn}

In \S\ref{sec-modifiedFenn} it is explained how to add 3-cells to a presentation in order to simplify diagrammatics in the presence of involutions and symmetries.

\begin{defn} Let $\PC = (\SC,\QC \cup \BC, \ZC \cup \MC)$ be the following 3-presentation, extending the usual presentation of $W$, which we call the \emph{Coxeter 3-presentation}. The
Zamolodzhikov 3-cells $\ZC$ are the same as in the previous chapter. The ``diagram-simplifying" 3-cells $\MC$ consist of: \begin{itemize} \item one 3-cell $z_s$ for each generating
involution $s \in \SC$, as in \eqref{RP2notcancel}; \item one 3-cell for each braid relation $(st)^m = 1$ accounting for rotational symmetry; \item one 3-cell for each braid relation $(st)^m = (ts)^m = 1$, accounting for flip symmetry. \end{itemize} The 3-cells accounting for rotational and flip symmetry were described in \S\ref{subsec-rotations}. \end{defn}

\begin{remark} As discussed in \S\ref{subsec-rotations}, the braid relation does not admit rotational or flip symmetry until one accounts for the fact that the generators are involutions.
\end{remark}

The following proposition is obvious from the definitions.

\begin{prop} $\tWC_\diag$ is isomorphic (not just equivalent) to $\Pi(\tX_\PC)^{\unor}_{\le 2}$ as $2$-categories. \end{prop}
	
Now we restate and prove one of our main theorems from the introduction.

\begin{thm} \label{thm:coxequiv} The obvious functor $\WC^\diag \to \Omega W$ is an equivalence of categories. \end{thm}

\begin{proof} It is enough to prove that $\pi_2(\tX_\PC)=0$. This follows from the lemma below. \end{proof}

\begin{lemma} By removing redundant 3-cells from $\tX_\PC$, one obtains a space which will deformation retract to the 3-skeleton of $\tCox$. In other words, $\tX_\PC$ is homotopy
equivalent to $\tCox^3 \vee S^3 \vee \cdots \vee S^3$. \end{lemma}

Example \ref{pancakeexample} illustrates the basic idea of this proof.

\begin{proof} We prove this lemma in steps. At the $k$-th step, we construct a sub-complex $\tX^{(k)}$ of $\tX_\PC$ by choosing certain cells to include. The sub-complex $\tX^{(k)}$ is not
the $k$-skeleton, though it will contain all $k$-cells of $\tX_\PC$ when $k<3$. We show that $\tX^{(k)}$ deformation retracts to $\tCox^k$. In particular, up to homotopy equivalence, we
can construct $\tX^{(k+1)}$ by gluing higher cells to $\tCox^k$ instead of $\tX^{(k)}$. For $k=3$, the difference between $\tX^{(3)}$ and $\tX_\PC^3$ will consist entirely of redundant
3-cells. Both $\tCox$ and $\tX_\PC$ have the same 0-skeleton, so we begin with $\tX^{(0)} = \tX^0_\PC$.

Now consider a single $s \in \SC$, and its parabolic subgroup $W_s \subset W$. By gluing in the 1-cells, 2-cells, and 3-cells corresponding to the sub-presentation $(s,s^2,z_s)$, one
obtains a copy of $S^3$ for each coset of $W_s$ in $W$ (see \S\ref{subsec-involutions}). One of the 3-cells is redundant, and excising it one obtains a copy of $D^3$ for each coset.
Each $D^3$ will deformation retract to a single edge between the two 0-cells, which can be thought of as an $s$-colored edge in $\tCox^1$. Thus if we take $\tX^{(0)}$ and add both 1-cells,
both 2-cells, and one 3-cell of $(s,s^2,z_s)$ for each $s \in \SC$, we obtain a space $\tX^{(1)}$ which deformation retracts to $\tCox^1$.

Now consider a single pair $s,t \in \SC$ with $m=m_{s,t} < \infty$, and its parabolic subgroup $W_{s,t} \subset W$. Let $b_{s,t}$ denote the braid relation inside $\BC$, $r_{s,t}$ denote
the rotation 3-cell inside $\MC$, and $f_{s,t}$ denote the flip 3-cell inside $\MC$. Each coset of $W_{s,t}$ in $W$ corresponds to a hollow $2m$-gon in $\tCox^1$ (or something which
deformation retracts to a hollow $2m$-gon in $\tX^{(1)}$). Gluing in the 2-cells corresponding to $b_{s,t}$, each coset will look like $2m$ disks, each glued along their boundary to a
common $S^1$. One can visualize this as an amalgamation of $2m-1$ copies of $S^2$, where the southern hemisphere of the $i$-th copy is identified with the northern hemisphere of the
$i+1$-st copy. There are a total of $4m$ 3-cells corresponding to $r_{s,t}$ and $f_{s,t}$ ($2m$ of each), each of which gives a cobordism between two different disks. One can choose $2m-1$
such 3-cells to fill in the $2m-1$ copies of $S^2$, yielding a cell complex structure on $D^3$. The remaining $(2m+1)$ 3-cells are all redundant, and we do not include them in $\tX^{(2)}$.
This copy of $D^3$ for each coset will deformation retract to a single solid $2m$-gon, which is a 2-cell in $\tCox^2$ corresponding to $W_{s,t}$. Thus if we take $\tX^{(1)}$ and add all
$2m$ 2-cells and $2m-1$ 3-cells of $(b_{s,t}, \{r_{s,t}, f_{s,t}\})$ for each $s,t \in \SC$ with $m_{s,t} < \infty$, we obtain a space $\tX^{(2)}$ which deformation retracts to $\tCox^2$.

Now consider a single triple $s,t,u \in \SC$ whose parabolic subgroup $W_{s,t,u}$ has finite size $n$. Let $Z_{s,t,u}$ denote the Zamolodzhikov 3-cell in $\ZC$. Each coset of $W_{s,t,u}$
gives a subspace of $\tCox^2$ which is a particular cell structure for $S^2$. In $\tX_\PC$, $Z_{s,t,u}$ corresponds to $n$ 3-cells (for each coset), each of which turns $S^2$ into $D^3$.
Clearly only one such 3-cell is necessary, after which the remaining ones are redundant. This single 3-cell corresponds precisely to the 3-cell in $\tCox^3$ for that coset. Thus if we take
$\tX^{(2)}$ and add a single 3-cell of the form $Z_{s,t,u}$ for each coset of $W_{s,t,u}$, we obtain the desired space $\tX^{(3)}$ which deformation retracts to $\tCox^3$. \end{proof}

\begin{remark} It is not unreasonable to expect a purely diagrammatic proof of Theorem \ref{thm:coxequiv}, and certainly this can be achieved in special cases. However, the difficulty in
finding this proof was what led the authors to this topological detour. \end{remark}


\bibliographystyle{plain}
\bibliography{everyone}{}

\end{document}